\documentclass[12pt,a4paper]{article}

\usepackage{amsmath}
\usepackage{amssymb}
\usepackage{amsthm}
\usepackage{amscd}
\usepackage{eepic}
\usepackage{graphicx}

\newtheorem{theorem}{Theorem}[section]
\newtheorem{lemma}[theorem]{Lemma}

\newtheorem{proposition}[theorem]{Proposition}
\newtheorem{remark}[theorem]{Remark}

\theoremstyle{definition}

\newtheorem{example}[theorem]{Example}

\makeatletter
\def\eqnarray{%
  \stepcounter{equation}%
  \let\@currentlabel=\theequation
  \global\@eqnswtrue
  \global\@eqcnt\z@
  \tabskip\@centering
  \let\\=\@eqncr
  $$\halign to \displaywidth\bgroup\@eqnsel\hskip\@centering
  $\displaystyle\tabskip\z@{##}$&\global\@eqcnt\@ne
  \hfil$\displaystyle{{}##{}}$\hfil
  &\global\@eqcnt\tw@$\displaystyle\tabskip\z@{##}$\hfil
  \tabskip\@centering&\llap{##}\tabskip\z@\cr}
\makeatother

\newcommand{\DIS}{\displaystyle}

\def\xv{\vec{x}}
\def\yv{\vec{y}}
\def\L{\mbox{\rm \boldmath{$L$}}}
\def\D{\mbox{\rm \boldmath{$D$}}}
\def\DIS{\displaystyle}

\def\cd{\cdots}
\def\ot{\otimes}

\def\bp{{\bf p}}

\def\veps{\varepsilon}
\def\vphi{\varphi}

\def\Z{{\mathbb Z}}
\def\et{\tilde{e}}
\def\ft{\tilde{f}}
\def\ve{\varepsilon}
\def\t{\mbox{\rm \boldmath{$t$}}}
\def\bp#1{\mbox{\mathversion{bold}$#1$}}

\title{The $A^{(1)}_M$ automata related to crystals of symmetric tensors}
\author{
G. Hatayama\thanks{
Institute of Physics, Graduate School of Arts and Sciences,
University of Tokyo, Komaba, Tokyo 153-8902, Japan},
K.  Hikami\thanks{
Department of Physics, Graduate School of Science,
University of Tokyo, Hongo, Tokyo 113-0033, Japan},
R. Inoue,$\hspace{-1.2mm}^\dagger$\\
A.  Kuniba$\hspace{-0.6mm}^*$,
T.  Takagi\thanks{
Department of Mathematics and Physics, National Defense Academy,
Yokosuka 239-8686, Japan},
and T. Tokihiro\thanks{
Graduate School of Mathematical Sciences, University of Tokyo, Komaba,
Tokyo 153-8914, Japan}}

\date{}

\begin{document}

\maketitle

\begin{abstract}
A soliton cellular automaton associated with crystals of symmetric tensor
representations of the quantum affine algebra
$U'_q(A^{(1)}_M)$ is introduced.
It is a crystal theoretic formulation of the
generalized box-ball system
in which capacities of boxes and carriers are
arbitrary and inhomogeneous.
Scattering matrices of two solitons
coincide with the combinatorial $R$
matrices of $U'_q(A^{(1)}_{M-1})$.
A piecewise linear evolution equation of the automaton is
identified with an ultradiscrete limit of the nonautonomous
discrete KP equation.
A class of $N$ soliton solutions is obtained through the
ultradiscretization of soliton solutions of the latter.
\end{abstract}

\section{Introduction}\label{sec:intro}

The box-ball system invented by
Takahashi and Satsuma  \cite{TS} is an important example of
soliton cellular automata.
It is a discrete dynamical system in which finitely many
balls move along the one dimensional array of
boxes under a certain rule.
Its integrability has been proved in \cite{TTMS}
by making a connection to the difference analogue of
the Lotka-Volterra equation \cite{HT}
through the limiting procedure called
{\em ultradiscrertization}.

By now the original box-ball system has been generalized into
several directions.
First, one can introduce the balls
distinguished by the index from the set $\{1,2,\ldots, M\}$.
Second, one lets the box at site $n$ accommodate up to $\theta_n$ balls,
where the capacity $\theta_n$ may depend on $n$.
Third, one can introduce a {\em carrier} with capacity $\kappa_t$
to redefine the time evolution at time $t$.
The carrier comes from the left and proceeds to the right,
picking up the balls in a box and dropping them into another
under a certain rule.
While it goes through the array of boxes, the successive
loading-unloading process induces the motion of balls over the boxes hence
the time evolution of the system.
These generalizations of the
Takahashi-Satsuma box-ball system are characterized by
the parameters $(M, \theta_n,  \kappa_t)$.
($n, t \in {\mathbb Z}$ play the role of space and time coordinates as in
the figure in Section \ref{subsec:cellauto}.)
The original one \cite{TS}
corresponds to the choice
$(M,\forall \theta_n, \forall \kappa_t) = (1, 1, \infty)$.
The case $(M, \forall \theta_n=1, \forall \kappa_t=\infty)$
was introduced in \cite{T} and studied in \cite{TNS}.
Similarly, the cases
$(M=1,\forall \theta_n=\theta, \forall \kappa = \kappa)$ with $\kappa > \theta$
and
$(M,\theta_n, \forall \kappa_t = \infty)$
were treated in \cite{TM} and \cite{TTM}, respectively.
These works have been done mainly from the viewpoint of the ultradiscretization.

The purpose of this paper is
to study the general $(M, \theta_n, \kappa_t)$ case.
In Section \ref{sec:ac} we formulate the corresponding
generalization of the box-ball system in terms of the crystal theory \cite{K,KMN1,KMN2}.
The latter is a representation theory of quantum groups at $q=0$.
The unexpected link between the crystals and the box-ball systems
has also been exploited in \cite{HKT,FOY} through a
crystal theoretic interpretation of the $L$-operator
approach \cite{HIK}.
The idea is to regard the box-ball system as a
solvable vertex model \cite{B} at $q=0$
under a  `ferromagnetic' boundary condition.
More concretely, the box-ball system corresponding to the data
$(M, \theta_n, \kappa_t)$ is naturally related to the
$U'_q(A^{(1)}_M)$ vertex model at $q=0$ whose inhomogeneity
in the quantum and auxiliary spaces is parametrized by
$\theta_n$'s and $\kappa_t$'s, respectively.

Let $B_l$ be the classical crystal of $U_q'(A^{(1)}_M)$
in the sense of \cite{KMN1} corresponding to the $l$-fold
symmetric tensor representation of $U_q(A_M)$.
Then the array of boxes and the ball configurations
are identified with
the elements from $\cdots \ot B_{\theta_{n}} \ot B_{\theta_{n+1}} \ot \cdots$.
The time evolution by the carrier with capacity $\kappa_t$
is realized as the action of the $q=0$ row transfer matrix
acting on $\cdots \ot B_{\theta_{n}} \ot B_{\theta_{n+1}} \ot \cdots$
with the  auxiliary space corresponding to $B_{\kappa_t}$.
We call the resulting dynamical system the $A^{(1)}_M$ automaton.
It is the most general one
in the $A^{(1)}_M$ case as far as the crystals for symmetric tensors are concerned.
For generalizations to other root systems,
see \cite{HI} for a supersymmetric one and
\cite{HKT} for the non exceptional series other than $A^{(1)}_M$.

In Section \ref{sec:proof} we introduce solitons and study the 2 soliton
scattering.
As in \cite{HKT, FOY} we label the solitons in terms of the elements of the
$U'_q(A^{(1)}_{M-1})$-crystal $B_l$, where $l$ plays the role of the amplitude of
a soliton.
In the collisions of two solitons associated with $B_l$ and $B_k$,
the scattering matrix is shown to coincide with the
combinatorial $R$ matrix giving the isomorphism
$B_l \ot B_k \simeq B_k \ot B_l$ of the $U'_q(A^{(1)}_{M-1})$-crystals.
These features are essentially the same with the
$\forall \theta_n = 1$ case \cite{TNS,HKT,FOY}.
A new aspect here is that depending on the
amplitudes $l, k$ and the parameters $\theta_n, \kappa_t$,
smaller soliton can overtake the larger one.
This is most transparently understood by viewing the
scattering from the cross channel.
By interchanging $\theta_n$ and $\kappa_t$, one can swap the
role of the space and time  and thereby
the boxes and carriers.
Then the curious scattering mentioned above reduces to the
`usual' one in the cross channel where the larger soliton overtakes
the smaller one.
In Section \ref{subsec:conserve} we also give a brief sketch of
the conserved quantities of our automaton following  \cite{FOY}.

In Section \ref{sec:toki} we set up piecewise linear
equation for the relevant combinatorial $R$ matrix \cite{NY} and
the resulting evolution equation for the $A^{(1)}_M$ automaton.
Extending the earlier result \cite{TTM},
we identify the evolution equation with an
ultradiscrete limit of the nonautonomous discrete
Kadomtsev-Petviashivili (ndKP) equation.
A class of $N$ soliton solutions is obtained through an
ultradiscretization of the $\tau$ functions.
As in the previous case \cite{TTM} one needs to make a
fine adjustment of the fermion momenta entering the
vacuum expectation value expression of the $\tau$ functions.
Each soliton in the automaton is obtained by
letting $M$ solitons in the ndKP
merge together in the ultradiscrete limit.

Section \ref{sec:summary} is a summary.
Appendices \ref{app:a} and \ref{app:b} contain the
details of the proofs of Proposition \ref{prop:toki}
and Theorem \ref{prop:N-soliton}, respectively.

\section{Automata from crystals}\label{sec:ac}

\subsection{$U'_q(A^{(1)}_M)$-crystals}
Let $B_k$ be the classical crystal of $U'_q(A^{(1)}_{M})$
corresponding to the $k$-fold symmetric tensor representation.
As a set it consists of the single row semistandard tableaux of length $k$
on letters $\{1,2,\ldots, M+1\}$:

\begin{equation*}
B_k = \{\fbox{$m_1 \cd m_k$}  \mid m_i \in \{1,\ldots, M+1\},
m_1 \le \cd \le m_k \},
\end{equation*}
where we have omitted the $k-1$ vertical lines separating the entries.
We also represent the elements by the multiplicities of their contents.
Namely, $b=\fbox{$m_1 \cd m_k$} \in B_k$ is also denoted by
$b=(x_1,x_2,\cd ,x_{M+1})$ with $x_i = \mbox{\#} \{ l \mid m_l = i\}$.

Denote the Kashiwara operators of $B_k$ by
$\ft_i$ and $\et_i$ for $i = 0,1,\ldots, M$.
The actions of $\et_i,\ft_i$ on $B_k$ are defined as
follows:
for $b=(x_1,x_2,\cd ,x_{M+1}) \in B_k$,
\begin{equation}\label{eq:action}
\left\{
	\begin{array}{ll}
	\et_0 b & =(x_1-1,x_2,\cd ,x_{M+1}+1),\\
	\ft_0 b & =(x_1+1,x_2,\cd ,x_{M+1}-1),\\
	\et_i b & =(x_1,\cd ,x_i+1,x_{i+1}-1,\cd ,x_{M+1})
	\quad \mbox{for}\, i=1, \cd ,M,\\
	\ft_i b & =(x_1,\cd ,x_i-1,x_{i+1}+1,\cd ,x_{M+1})
	\quad \mbox{for}\, i=1, \cd ,M.
	\end{array}
\right.
\end{equation}
In the above, the right hand sides are to be understood as $0$
if they are not in $B_k$.
A crystal can be regarded as a colored oriented graph called a ``crystal graph"
by defining
\[
b\stackrel{i}{\longrightarrow}b'\quad\Longleftrightarrow\quad \ft_{i}b=b'.
\]
Thus for example
$B_1 = \{\fbox{1}, \ldots, \fbox{$M\!+\!1$}\}$
has the crystal graph:\\
\setlength{\unitlength}{0.35mm}
\begin{picture}(300,80)(0,-10)
\put(10,10){
\bezier{500}(10,20)(147.5,60)(285,20)
\put(10,20){\vector(-3,-1){1}}
\multiput(20,10)(50,0){3}{\vector(1,0){30}}
\put(5,5){\makebox(10,10){\fbox{1}}}
\put(55,5){\makebox(10,10){\fbox{2}}}
\put(105,5){\makebox(10,10){\fbox{3}}}
\put(32,10){\makebox(10,10){1}}
\put(81,10){\makebox(10,10){2}}
\put(130,10){\makebox(10,10){3}}
\multiput(157,10)(5,0){6}{\circle*{1}}
\multiput(190,10)(55,0){2}{\vector(1,0){30}}
\put(225,5){\makebox(15,10){\fbox{$M$}}}
\put(289,5){\makebox(10,10){\fbox{$M\!+\!1$}.}}
\put(199,10){\makebox(15,10){$M\!-\!1$}}
\put(255,10){\makebox(10,10){$M$}}
\put(147.5,40){\makebox(10,10){0}}
}
\end{picture}

\begin{figure}[htbp]
\begin{center}
\includegraphics[width=.6\linewidth]{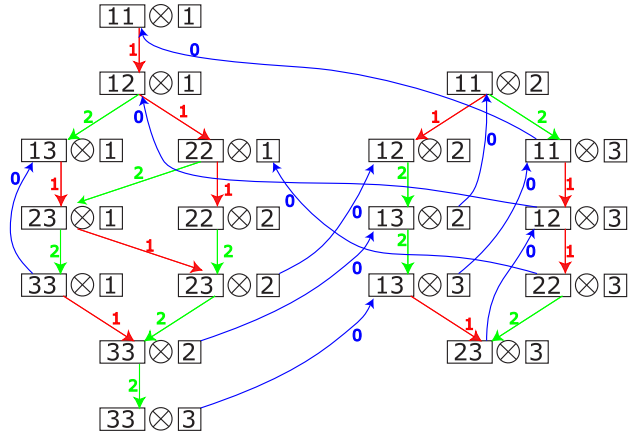}
\end{center}
\caption{crystal graph of $U'_q(A^{(1)}_2)$-crystal $B_2 \ot B_1$}
\label{fig:A2B2B1}
\end{figure}

\begin{figure}[htbp]
\begin{center}
\includegraphics[width=.6\linewidth]{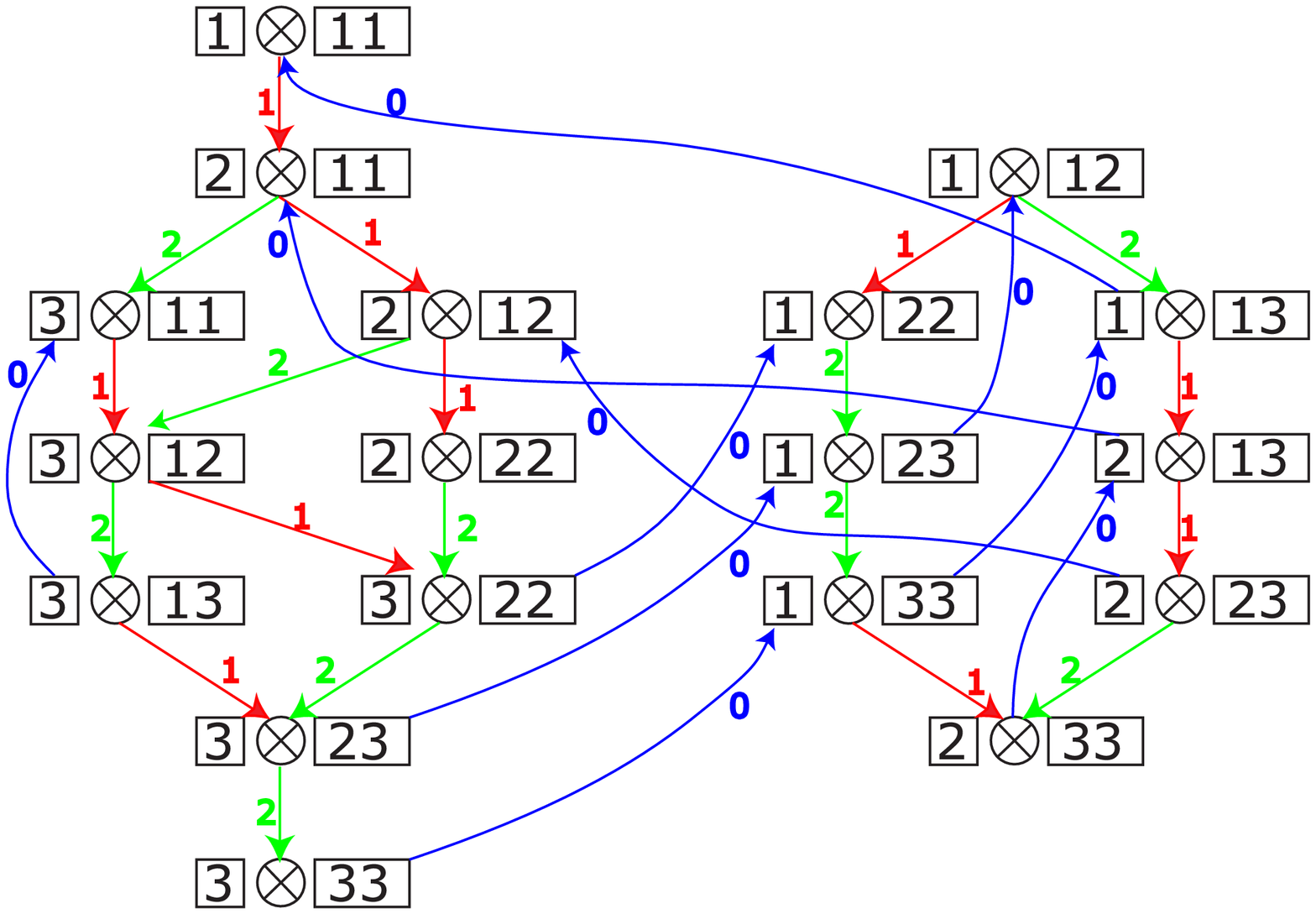}
\end{center}
\caption{crystal graph of $U'_q(A^{(1)}_2)$-crystal $B_1 \ot B_2$}
\label{fig:A2B1B2}
\end{figure}

Setting $\veps_i(b) = \max_l\{\et_i^l b \neq 0 \mid l \ge 0 \}$ and
$\vphi_i(b) = \max_l\{\ft_i^l b \neq 0 \mid l \ge 0 \}$
for $b \in B_k$, one has
\begin{alignat*}{2}
\veps_0(b) &= x_1, \quad & \veps_i(b) &= x_{i+1}  \quad \mbox{for}\, i=1,
\cd ,M,\\
\vphi_0(b) &= x_{M+1}, \quad & \vphi_i(b) &= x_i  \quad \mbox{for}\, i=1,
\cd ,M.
\end{alignat*}
This data is necessary when we treat tensor products of the crystals.
For two crystals $B$ and $B'$, the tensor product $B\ot B'$
is defined. As a set,
\[
B\ot B'=\{b_1\ot b_2\mid b_1\in B,b_2\in B'\}.
\]
The actions of $\et_{i}$ and $\ft_{i}$ are defined by
\begin{eqnarray}
\et_{i}(b_1\ot b_2)&=&\left\{
\begin{array}{ll}
\et_{i}b_1\ot b_2&\mbox{ if }\vphi_i(b_1)\ge\veps_i(b_2)\\
b_1\ot \et_{i}b_2&\mbox{ if }\vphi_i(b_1) < \veps_i(b_2),
\end{array}\right. \label{eq:ot-e}\\
\ft_{i}(b_1\ot b_2)&=&\left\{
\begin{array}{ll}
\ft_{i}b_1\ot b_2&\mbox{ if }\vphi_i(b_1) > \veps_i(b_2)\\
b_1\ot \ft_{i}b_2&\mbox{ if }\vphi_i(b_1)\le\veps_i(b_2).
\end{array}\right. \label{eq:ot-f}
\end{eqnarray}
Here $0\ot b$ and $b\ot0$ are understood to be $0$.
For two crystals $B$ and $B'$,
the tensor products $B' \otimes B$ and $B \otimes B'$
constructed as above
are again crystals which are canonically isomorphic.
The isomorphism $R \, : \, B' \otimes B
\stackrel{\sim}{\rightarrow} B\otimes B'$ is called the
combinatorial $R$ matrix \cite{KMN1,NY}.
By the definition $R$ commutes with $\ft_i, \et_i$ for any
$i = 0,1, \ldots, M$.
(More precisely one introduces affine crystals and the associated
energy function, but in this paper we shall exclusively
treat classical crystals and concern the energy function only
in connection with the conserved quantities in Section \ref{subsec:conserve}.)
\begin{example}\label{ex:crystalgraph}
Figure \ref{fig:A2B2B1} and Figure \ref{fig:A2B1B2} are
the crystal graphs of
$U'_q(A^{(1)}_2)$-crystals $B_2 \ot B_1$ and
$B_1 \ot B_2$, respectively.
\end{example}
\begin{example}\label{ex:R}
Let $B'=B_2,\,B=B_1$ of $U'_q(A^{(1)}_2)$-crystals.
\begin{eqnarray*}
\makebox{ (i) }  && R: \fbox{$13$} \ot \fbox{$2$} \mapsto \fbox{$1$} \ot
\fbox{$23$},\\
\makebox{(ii) } && R: \fbox{$23$} \ot \fbox{$2$} \mapsto \fbox{$3$} \ot
\fbox{$22$}.
\end{eqnarray*}
These are obtained by comparing the crystal graphs
in Example \ref{ex:crystalgraph}.
\end{example}

%

%
We write the highest weight element in $B_k$ with respect to
$U_q(A_M)$ as $u_k$:
\begin{equation}\label{eq:udef}
u_k = \overbrace{\fbox{$1 1 \ldots 1$}}^k
=(k,0,\cd ,0) \in B_k.
\end{equation}
%

\subsection{Isomorphism}
Here we give an explicit procedure to obtain
the isomorphism $R : B_k \ot B_l \rightarrow B_l \ot B_k$
without drawing the whole crystal graphs of $B_k \ot B_l$ and
$B_l \ot B_k$.

Let $b_1 \ot b_2$ be an element in $B_k \ot B_l$
such as $b_1=(x_1,\ldots,x_{M+1})$ and $b_2=(y_1,\ldots,y_{M+1})$.
We represent $b_1 \ot b_2$ by the two column diagram.
Each column has $M+1$ rows, enumerated as 1 to $M+1$
{}from the top to the bottom.
We put $x_i$ (resp. $y_i$) dots $\bullet$ in the $i$-th row
of the left (resp. right) column.

\setlength{\unitlength}{9mm}
\begin{picture}(5,4.5)(-5,0)
\put(-2.5,1.5){\makebox(2,1){$b_1 \ot b_2=$}}
\put(0,0){\line(0,1){4}}
\put(2,0){\line(0,1){4}}
\put(0,0){\line(1,0){2}}
\put(0,0){\makebox(2,1){$\scriptstyle{\underbrace{\bullet \cdots \bullet}
_{x_{M+1}}}$}}
\put(0,1){\line(1,0){2}}
\put(0,1){\makebox(2,1){$\vdots$}}
\put(0,2){\line(1,0){2}}
\put(0,2){\makebox(2,1){$\scriptstyle{\underbrace{\bullet \cdots \bullet}
_{x_{2}}}$}}
\put(0,3){\line(1,0){2}}
\put(0,3){\makebox(2,1){$\scriptstyle{\underbrace{\bullet \cdots \bullet}
_{x_{1}}}$}}
\put(0,4){\line(1,0){2}}
\put(3,0){\line(0,1){4}}
\put(5,0){\line(0,1){4}}
\put(3,0){\line(1,0){2}}
\put(3,0){\makebox(2,1){$\scriptstyle{\underbrace{\bullet \cdots \bullet}
_{y_{M+1}}}$}}
\put(3,1){\line(1,0){2}}
\put(3,1){\makebox(2,1){$\vdots$}}
\put(3,2){\line(1,0){2}}
\put(3,2){\makebox(2,1){$\scriptstyle{\underbrace{\bullet \cdots \bullet}
_{y_{2}}}$}}
\put(3,3){\line(1,0){2}}
\put(3,3){\makebox(2,1){$\scriptstyle{\underbrace{\bullet \cdots \bullet}
_{y_{1}}}$}}
\put(3,4){\line(1,0){2}}
\end{picture}

\begin{proposition}\label{prop:H-rule}
The rule to obtain the
isomorphism $R$ is as follows.
\begin{itemize}
\item[(1)]
Assume $k \geq l$ (resp.~$k \leq l$).
Pick any dot, say $\bullet_a$, in the right (resp.~left)
column and find its partner $\bullet_a'$
in the left (resp.~right) column.
The $\bullet_a'$ is chosen
{}from the dots which are in the lowest (resp.~highest) row among all dots
whose positions are higher (resp.~lower) than that of $\bullet_a$.
If there is no such dot, we return to the bottom (resp.~top) and
the partner $\bullet_a'$ is chosen {}from the dots
in the lowest (resp.~highest) row among all dots.
Connect $\bullet_a$ and $\bullet_a'$ by a line.
We call the lines in the latter case winding and in the former case unwinding.
\item[(2)]
Repeat the procedure (1) for the remaining unconnected dots
$(l-1)$-times (resp.~$(k-1)$-times).
\item[(3)]
The isomorphism $R$ is obtained by sliding
the remaining $(k-l)$ (resp.~$(l-k)$) unpaired dots in the
left (resp.~right) column to the right (resp.~left).
\end{itemize}
\end{proposition}
\par\noindent
The $R$ obtained by this rule has the correct property
as the isomorphism.
This fact has been proved in section 3 of \cite{NY}.
We will write the relation $R: u \ot v \mapsto v' \ot u'$
also as $u \ot v \simeq v' \ot u'$.
Obviously one has
\begin{equation}\label{eq:hh}
u_k \ot u_l \simeq u_l \ot u_k
\end{equation}
for the element (\ref{eq:udef}).
\begin{example}\label{ex:NY}
Let $M=2,\,k=2,\,l=1$.
Example \ref{ex:R} (i),(ii) are obtained
by the following diagrams:
\begin{center}
\raisebox{8mm}{(i)}
\setlength{\unitlength}{3mm}
\begin{picture}(6,6)(0,0)
\multiput(0,0)(0,2){4}{\line(1,0){2}}
\multiput(4,0)(0,2){4}{\line(1,0){2}}
\multiput(0,0)(2,0){4}{\line(0,1){6}}
\put(1,1){\circle*{.6}}
\put(1,5){\circle*{.6}}
\put(5,3){\circle*{.6}}
\path(5,3)(1,5)
\end{picture}
\raisebox{8mm}{,\hspace{10pt}(ii)}
\setlength{\unitlength}{3mm}
\begin{picture}(6,6)(0,0)
\multiput(0,0)(0,2){4}{\line(1,0){2}}
\multiput(4,0)(0,2){4}{\line(1,0){2}}
\multiput(0,0)(2,0){4}{\line(0,1){6}}
\put(1,1){\circle*{.6}}
\put(1,3){\circle*{.6}}
\put(5,3){\circle*{.6}}
\path(5,3)(1,1)
\end{picture}
\raisebox{8mm}{.}
\end{center}
The line in (i) is unwinding   and that in (ii) is winding.
\end{example}

Suppose $b\ot b'\in B\ot B'$ is mapped to $\tilde{b'}\ot \tilde{b}
\in B'\ot B$ under the isomorphism
$B \ot B' \simeq B' \ot B$
of $U_q'(A^{(1)}_M)$-crystals.
A $\Z$-valued function
$H$ on $B\ot B'$ is called an {\em energy function} if for any $i$
and $b\ot b'\in B\ot B'$ such that $\et_i(b\ot b')\neq0$,
 it satisfies
\begin{eqnarray}
H(\et_i(b\ot b'))&=H(b\ot b')+1
&\mbox{ if }i=0,\vphi_0(b)\geq\veps_0(b'),\vphi_0(\tilde{b'})\geq\veps_0(\tilde{b}),
\nonumber\\
&=H(b\ot b')-1
&\mbox{ if }i=0,\vphi_0(b)<\veps_0(b'),\vphi_0(\tilde{b'})<\veps_0(\tilde{b}),
\nonumber\\
&\hspace{-6mm}=H(b\ot b')&\mbox{ otherwise}.\label{eq:e-func}
\end{eqnarray}
When we want to emphasize $B\ot B'$, we write $H_{BB'}$ for $H$.
This definition of the energy function
is due to (3.~4.~e) of \cite{NY}, that is a
generalization of the definition for $B=B'$ case
in \cite{KMN1}.
The energy function is unique up to additive constant, since
$B\ot B'$ is connected. By definition, $H_{BB'}(b\ot b')
=H_{B'B}(\tilde{b'}\ot \tilde{b})$.
Throughout this paper we normalize it as
\begin{equation}\label{eq:Hnormalize}
H_{B_lB_k}(u_l \ot u_k) =0,
\end{equation}
irrespective of $l < k$ or $l \ge k$.
Then it is the result of \cite{NY} that
the energy function is $(-1)$ times the number of
unwinding lines in the sense  of Example \ref{ex:NY}.

With a successive application of $R$'s, one interchanges
the order of tensor product pairwise and obtains the
isomorphism of $B_{k_1}\ot \cd \ot B_{k_n}$ and
$B_{k_{P_1}}\ot \cd \ot B_{k_{P_n}}$ for any permutation $P$.
The compatibility of this construction is guaranteed by the
Yang-Baxter equation obeyed by $R$.
The following assertion follows easily from Proposition \ref{prop:H-rule}.
\begin{proposition}\label{prop:ferro}
Let $k_1, k_2, \ldots\in {\mathbb Z}_{\ge 1}$
be any sequence.
Suppose $b \ot u_{k_1} \ot \cd \ot u_{k_n} \simeq c_1 \ot \cd \ot c_n \ot b'$
is valid for some $b'$ and $c_i$'s under the isomorphism
$B_l \ot B_{k_1} \ot \cd \ot B_{k_n} \simeq
B_{k_1} \ot \cd \ot B_{k_n} \ot B_l$.
For any $b \in B_l$, there exists $n_0$ such that $b'= u_l$ for all $n \ge n_0$.
\end{proposition}
This property will be needed in
constructing the automaton in Section \ref{subsec:cellauto}.
\subsection{Automaton}
\label{subsec:cellauto}
Let
$\cd ,\theta_{-1},\theta_0,\theta_1,\cd$ and
$\cd ,\kappa_{-1},\kappa_0,\kappa_1,\cd$ be
two sequences of positive integers.
Denote the former indices by $n$, and the latter
indices by $t$.
Consider the 2D-lattice with $n$- and $t$- directions,
\begin{eqnarray*}
&&\mbox{\rm $n$-direction} \qquad \cdots \otimes B_{\theta_{n-1}} \otimes
B_{\theta_n} \otimes B_{\theta_{n+1}} \otimes \cdots \\
&&\mbox{\rm $t$-direction} \qquad \cdots \otimes B_{\kappa_{t-1}} \otimes
B_{\kappa_t} \otimes B_{\kappa_{t+1}} \otimes \cdots .
\end{eqnarray*}
In terms of the box-ball systems, $\theta_n$ is the
capacity of the $n$-th boxCand $\kappa_t$ is the capacity of the $t$-th
carrier.

Draw $t$- constant lines horizontally, and $n$- constant lines
vertically.
Number the former downward, and the latter to the right.
At any horizontal or vertical
line segment of the lattice, we inscribe an element of the crystals
in the following way.
At the point labeled by $(t,n)$,
we put $b^t_n \in B_{\theta_n}$ on the upper line segment and
$v^t_n \in B_{\kappa_t}$ on the left line segment.
Thus we have
$b^{t+1}_n \in B_{\theta_n}$ on the lower line segment and
$v^t_{n+1} \in B_{\kappa_t}$ on the right line segment.

\setlength{\unitlength}{1mm}
\begin{picture}(100,50)(-30,0)
\put(0,10){\line(1,0){40}}
\put(0,20){\line(1,0){40}}
\put(0,30){\line(1,0){40}}
\put(10,0){\line(0,1){40}}
\put(20,0){\line(0,1){40}}
\put(30,0){\line(0,1){40}}
\put(-6,7.5){\makebox(5,5){$\scriptstyle{t+1}$}}
\put(-6,17.5){\makebox(5,5){$\scriptstyle{t}$}}
\put(-6,27.5){\makebox(5,5){$\scriptstyle{t-1}$}}
\put(7.5,40){\makebox(5,5){$\scriptstyle{n-1}$}}
\put(17.5,40){\makebox(5,5){$\scriptstyle{n}$}}
\put(27.5,40){\makebox(5,5){$\scriptstyle{n+1}$}}
\put(14,12.2){\makebox(5,5){$b^{t+1}_n$}}
\put(15,22.8){\makebox(5,5){$b^{t}_n$}}
\put(10.8,20.3){\makebox(5,5){$v^t_{n}$}}
\put(22.7,20.3){\makebox(5,5){$v^t_{n+1}$}}
\end{picture}
\par\noindent
We impose the condition that they are related by the combinatorial $R$ matrix,
\begin{equation}\label{eq:vbbv}
R \, : \, v^t_{n} \ot b^{t}_n
\stackrel{\sim}{\mapsto} b^{t+1}_n \ot v^t_{n+1}.
\end{equation}
In the following sections,
we  consider the time evolution of the system downward.
In view of Proposition \ref{prop:ferro}
we can and will exclusively consider the case where
for any $t$, $b^t_n \neq u_{\theta_n}$ only for finitely many $n$'s
and similarly for any $n$, $v^t_n \neq u_{\kappa_t}$ only for finitely many $t$'s.
Sometimes we ignore $v^t_n$'s and display
the time evolution of the system only with  the arrays
\begin{equation*}
\begin{array}{ccccccc}
\ldots & b^0_{-2} & b^0_{-1} & b^0_0 & b^0_1 & b^0_2 &\ldots \\
\ldots & b^1_{-2} & b^1_{-1} & b^1_0 & b^1_1 & b^1_2 &\ldots \\
\ldots & b^2_{-2} & b^2_{-1} & b^2_0 & b^2_1 & b^2_2 &\ldots .\\
\end{array}
\end{equation*}
In short, the evolution of the array $\{b_n^{t}\}$ to $\{b_n^{t+1}\}$
is determined by
\[
\begin{array}{rcl}
B_{\kappa_t} \ot ( \cdots \ot B_{\theta_{n}} \ot
B_{\theta_{n+1}} \ot \cdots ) &
\simeq &
( \cdots \ot B_{\theta_{n}} \ot
B_{\theta_{n+1}} \ot \cdots ) \ot B_{\kappa_t}
\\
u_{\kappa_t} \ot ( \cdots \ot b_{n}^{t} \ot
b_{n+1}^{t} \ot \cdots ) &
\simeq &
( \cdots \ot b_{n}^{t+1} \ot
b_{n+1}^{t+1} \ot \cdots ) \ot u_{\kappa_t} \\
\end{array}
\]
under the successive applications of the combinatorial $R$ matrices
$R: B_{\kappa_t} \ot B_{\theta_j} \xrightarrow{\sim}
B_{\theta_j} \ot B_{\kappa_t}$.

Setting
$\mbox{\mathversion{bold}$p$} = \cdots \ot b_{n}^{t} \ot
b_{n+1}^{t} \ot \cdots$, we denote the time evolution induced by $u_{\kappa_t}$
as above by $T_{\kappa_t}(\mbox{\mathversion{bold}$p$}) = \cdots \ot b_{n}^{t+1} \ot
b_{n+1}^{t+1} \ot \cdots$.
Obviously the time evolutions are invertible, and due to (\ref{eq:hh})
they are commutative,
\begin{equation}\label{eq:commute}
T_{\kappa} T_{\kappa'} = T_{\kappa'} T_{\kappa}.
\end{equation}
In the rest of the paper, the 2 dimensional lattice
on which the automaton is defined should be appropriately
understood either as large but finite
or formally infinite  depending on the situation.

The following observation will turn out useful in the sequel.
\begin{remark}\label{rem:pc}
Interchanging the role of `space' and `time',
one can view (\ref{eq:vbbv}) as the evolution of the
array $\cd \ot v^{t+1}_{n+1} \ot v^t_{n+1} \ot v^{t-1}_{n+1} \ot \cd$ to the left as
\begin{equation*}
T_{\theta_n}(\cd \ot v^{t+1}_{n+1} \ot v^t_{n+1} \ot v^{t-1}_{n+1} \ot \cd)
= \cd \ot v^{t+1}_{n}\ot v^t_{n} \ot v^{t-1}_{n} \ot \cd.
\end{equation*}
\end{remark}
\begin{example}\label{ex:tns}
Let $M=3, \forall \theta_n=1$ and $\forall \kappa_t = \infty$.\hfill
\begin{center}
$\cdots$111142113111111111111$\cdots$\\
$\cdots$111111421311111111111$\cdots$\\
$\cdots$111111114231111111111$\cdots$\\
$\cdots$111111111124311111111$\cdots$\\
$\cdots$111111111112143111111$\cdots$\\
$\cdots$111111111111211431111$\cdots$\\
\end{center}
where $i$ denotes \fbox{$i$}.
This is a typical 2 soliton scattering.
One can see that a soliton with amplitude $l$
moves to the right with velocity $l$ if separated sufficiently.
Hence the larger solitons overtake the smaller ones.
(See Section \ref{subsec:solitons} for the precise definition of the
solitons and their amplitude.)
\end{example}
\begin{example}\label{ex:tnsd}
Let $M=3, \forall \theta_n = 2$ and $\forall \kappa_t = 1$.

\setlength{\unitlength}{1mm}
\begin{picture}(135,95)(10,-5)
\put(5,0.5){11}
\put(18,0.5){11}
\put(31,0.5){11}
\put(44,0.5){11}
\put(57,0.5){11}
\put(70,0.5){11}
\put(83,0.5){34}
\put(96,0.5){11}
\put(109,0.5){11}
\put(122,0.5){11}
\put(135,0.5){13}

\multiput(2.5,7)(13,0){11}{\line(1,0){9}}
\multiput(7,4.5)(13,0){11}{\line(0,1){5}}
\put(-0.5,6){1}
\put(12.5,6){1}
\put(25.5,6){1}
\put(38.5,6){1}
\put(51.5,6){1}
\put(64.5,6){1}
\put(77.5,6){4}
\put(90.5,6){1}
\put(103.5,6){1}
\put(116.5,6){1}
\put(129.5,6){3}
\put(142.5,6){1}

\put(5,10.5){11}
\put(18,10.5){11}
\put(31,10.5){11}
\put(44,10.5){11}
\put(57,10.5){11}
\put(70,10.5){14}
\put(83,10.5){13}
\put(96,10.5){11}
\put(109,10.5){11}
\put(122,10.5){13}
\put(135,10.5){11}

\multiput(2.5,17)(13,0){11}{\line(1,0){9}}
\multiput(7,14.5)(13,0){11}{\line(0,1){5}}
\put(-0.5,16){1}
\put(12.5,16){1}
\put(25.5,16){1}
\put(38.5,16){1}
\put(51.5,16){1}
\put(64.5,16){1}
\put(77.5,16){3}
\put(90.5,16){1}
\put(103.5,16){1}
\put(116.5,16){3}
\put(129.5,16){1}
\put(142.5,16){1}

\put(5,20.5){11}
\put(18,20.5){11}
\put(31,20.5){11}
\put(44,20.5){11}
\put(57,20.5){11}
\put(70,20.5){34}
\put(83,20.5){11}
\put(96,20.5){11}
\put(109,20.5){13}
\put(122,20.5){11}
\put(135,20.5){11}

\multiput(2.5,27)(13,0){11}{\line(1,0){9}}
\multiput(7,24.5)(13,0){11}{\line(0,1){5}}
\put(-0.5,26){1}
\put(12.5,26){1}
\put(25.5,26){1}
\put(38.5,26){1}
\put(51.5,26){1}
\put(64.5,26){4}
\put(77.5,26){1}
\put(90.5,26){1}
\put(103.5,26){3}
\put(116.5,26){1}
\put(129.5,26){1}
\put(142.5,26){1}

\put(5,30.5){11}
\put(18,30.5){11}
\put(31,30.5){11}
\put(44,30.5){11}
\put(57,30.5){14}
\put(70,30.5){13}
\put(83,30.5){11}
\put(96,30.5){13}
\put(109,30.5){11}
\put(122,30.5){11}
\put(135,30.5){11}

\multiput(2.5,37)(13,0){11}{\line(1,0){9}}
\multiput(7,34.5)(13,0){11}{\line(0,1){5}}
\put(-0.5,36){1}
\put(12.5,36){1}
\put(25.5,36){1}
\put(38.5,36){1}
\put(51.5,36){1}
\put(64.5,36){3}
\put(77.5,36){1}
\put(90.5,36){3}
\put(103.5,36){1}
\put(116.5,36){1}
\put(129.5,36){1}
\put(142.5,36){1}

\put(5,40.5){11}
\put(18,40.5){11}
\put(31,40.5){11}
\put(44,40.5){11}
\put(57,40.5){34}
\put(70,40.5){11}
\put(83,40.5){13}
\put(96,40.5){11}
\put(109,40.5){11}
\put(122,40.5){11}
\put(135,40.5){11}

\multiput(2.5,47)(13,0){11}{\line(1,0){9}}
\multiput(7,44.5)(13,0){11}{\line(0,1){5}}
\put(-0.5,46){1}
\put(12.5,46){1}
\put(25.5,46){1}
\put(38.5,46){1}
\put(51.5,46){4}
\put(64.5,46){1}
\put(77.5,46){3}
\put(90.5,46){1}
\put(103.5,46){1}
\put(116.5,46){1}
\put(129.5,46){1}
\put(142.5,46){1}

\put(5,50.5){11}
\put(18,50.5){11}
\put(31,50.5){11}
\put(44,50.5){14}
\put(57,50.5){13}
\put(70,50.5){13}
\put(83,50.5){11}
\put(96,50.5){11}
\put(109,50.5){11}
\put(122,50.5){11}
\put(135,50.5){11}

\multiput(2.5,57)(13,0){11}{\line(1,0){9}}
\multiput(7,54.5)(13,0){11}{\line(0,1){5}}
\put(-0.5,56){1}
\put(12.5,56){1}
\put(25.5,56){1}
\put(38.5,56){4}
\put(51.5,56){1}
\put(64.5,56){3}
\put(77.5,56){1}
\put(90.5,56){1}
\put(103.5,56){1}
\put(116.5,56){1}
\put(129.5,56){1}
\put(142.5,56){1}

\put(5,60.5){11}
\put(18,60.5){11}
\put(31,60.5){14}
\put(44,60.5){11}
\put(57,60.5){33}
\put(70,60.5){11}
\put(83,60.5){11}
\put(96,60.5){11}
\put(109,60.5){11}
\put(122,60.5){11}
\put(135,60.5){11}

\multiput(2.5,67)(13,0){11}{\line(1,0){9}}
\multiput(7,64.5)(13,0){11}{\line(0,1){5}}
\put(-0.5,66){1}
\put(12.5,66){1}
\put(25.5,66){4}
\put(38.5,66){1}
\put(51.5,66){3}
\put(64.5,66){1}
\put(77.5,66){1}
\put(90.5,66){1}
\put(103.5,66){1}
\put(116.5,66){1}
\put(129.5,66){1}
\put(142.5,66){1}

\put(5,70.5){11}
\put(18,70.5){14}
\put(31,70.5){11}
\put(44,70.5){13}
\put(57,70.5){13}
\put(70,70.5){11}
\put(83,70.5){11}
\put(96,70.5){11}
\put(109,70.5){11}
\put(122,70.5){11}
\put(135,70.5){11}

\multiput(2.5,77)(13,0){11}{\line(1,0){9}}
\multiput(7,74.5)(13,0){11}{\line(0,1){5}}
\put(-0.5,76){1}
\put(12.5,76){4}
\put(25.5,76){1}
\put(38.5,76){1}
\put(51.5,76){3}
\put(64.5,76){1}
\put(77.5,76){1}
\put(90.5,76){1}
\put(103.5,76){1}
\put(116.5,76){1}
\put(129.5,76){1}
\put(142.5,76){1}

\put(5,80.5){14}
\put(18,80.5){11}
\put(31,80.5){11}
\put(44,80.5){33}
\put(57,80.5){11}
\put(70,80.5){11}
\put(83,80.5){11}
\put(96,80.5){11}
\put(109,80.5){11}
\put(122,80.5){11}
\put(135,80.5){11}
\end{picture}
Here $i$ and $ij$ denote \fbox{$i$} and \fbox{$i j$}, respectively.
We have depicted the both variables $\{b^t_n \}$ and $\{ v^t_n\}$.
This time  $14$
on the top left is the smaller soliton
and $33$ or $13 \ot 13$ is the larger soliton.
Thus in terms of the $\{b^t_n\}$ variable,
the smaller one overtakes the larger one as we go down the
figure ending with the solitons $34$ and $13$.
This is an opposite feature from the previous example.
However in the space-time interchanged picture (Remark \ref{rem:pc}),
it reduces to the situation similar to Example \ref{ex:tns}.
Namely, in terms of the $\{ v^t_n\}$ variable,
the larger soliton overtakes the small one as
$\ldots 43 \ldots 3 \ldots \rightarrow
\ldots 4 \ldots 33 \ldots$, as we trace the diagram from the
right to the left.
\end{example}

\subsection{Equivalence with box-ball systems}\label{subsec:equiva}

Our $A^{(1)}_M$ automaton can be viewed as a generalized box-ball system.
One interprets the letter 1 in the tableaux as an empty space and
the other letters $2 \le i \le M+1$ as the balls with index $M+2-i$.
The element $b^t_n$ signifies the balls contained in the $n$ th
box with capacity $\theta_n$ at time $t$.
Similarly $v^t_n$ stands for the carrier with capacity $\kappa_t$.
Then (\ref{eq:vbbv}) tells that through the loading-unloading process,
the box and the carrier change into
$b^{t+1}_n$ and $v^t_{n+1}$, respectively.
Sending the carrier through to the left, one has the
time evolution of the box-ball state
$\cd \ot b^t_n \ot b^t_{n+1} \ot \cd$ into
$\cd \ot b^{t+1}_n \ot b^{t+1}_{n+1} \ot \cd$.
For a concrete rule describing (\ref{eq:vbbv}) in terms of the box-ball
terminology, see {\em BBS scattering rule} in \cite{TNS}.
The relation (\ref{eq:vbbv}) will also be expressed as a piecewise linear
equation in Proposition \ref{prop:BBS-rule}.

When $\forall \kappa_t = \infty$ we claim that
the evolution of $\{b^t_n \}$ in
our $A^{(1)}_M$ automaton is equivalent
to the box-ball system studied in \cite{TTM} under the
above stated translation.
In the latter the one-dimensional array of boxes with  capacities
$\ldots, \theta_{n-1}, \theta_n, \theta_{n+1},\ldots$ accommodate the balls
with an index from the set $\{1, \ldots, M\}$.
The dynamics of the balls in each time step is governed by the rule \cite{TTM}:
\begin{enumerate}
\item Move every ball only once.
\item Move the leftmost ball with index 1 to the nearest right box with space.
\item Move the leftmost ball with index 1 among the rest to its nearest right box
with space.
\item Repeat this procedure until all of the balls with index 1 are moved.
\item Do the same procedure $2-4$ for the balls with index 2.
\item Repeat this procedure successively until all of the balls with index $M$ are moved.
\end{enumerate}
If the ball with some index is absent,
one just proceeds to those with the next index.
A box with space means the one that contains strictly
fewer balls than its capacity.
If a box contains more than one balls with the same index and they are not yet moved
at an instant during the procedure, one may pick any one of them
when looking for the leftmost one.
The equivalence to our automaton with
$\forall \kappa_t \rightarrow \infty$ is shown by the fact that the both
lead to the same evolution equation, which is
given {}from  Proposition \ref{prop:BBS-rule} under the said limit.

The above rule tells that
the time evolution $T_\infty$ in our automaton
admits the factorization:
\begin{equation}\label{eq:factor}
T_\infty = {\tilde T}_M \cdots {\tilde T}_2{\tilde T}_1,
\end{equation}
where ${\tilde T}_j$ moves the balls with index $j$ only, and we
identify the left hand side with the corresponding
operator acting on the box-ball systems.

For a later convenience we introduce the {\em canonical system}
following \cite{TNS}.
We keep assuming $\forall \kappa_t = \infty$
and stay in the description in terms of the box-ball system rather than
crystals until the end of this subsection.
Thus we identify $b \in B_\theta$ with the capacity $\theta$ box containing
the balls as specified before.
Suppose a state $\bp{p} = \cd \ot b_n \ot b_{n+1} \ot \cd$ contains
$J$ balls in total.
Then the action of ${\tilde T}_M \cdots {\tilde T}_2{\tilde T}_1$
consists of $J$ steps, each of which is to move a certain ball.
To a ball to be moved in the $j$ th step ($1 \le j \le J$),
we assign a {\em signature} $j$.
The assignment is unique up to the trivial freedom among the commonly
indexed balls within the same boxes.
Let $c(\bp{p})$ be the ball configurations obtained from $\bp{p}$ just by
regarding the signatures as  new indices.
It consists of the same array of the boxes  and
$J$ balls as before but with the new distinct index from $1$ to $J$.
One can still let $c(\bp{p})$ evolve under the previously stated rule $1-6$
by replacing $M$ therein  with $J$.
The resulting new box-ball system is called the canonical system.
By a  close inspection of the rule 1-6, it is not difficult to confirm
the commutativity:
\begin{equation}\label{eq:canonicalcom}
c({\tilde T}_M \cdots {\tilde T}_2{\tilde T}_1(\bp{p}))
={\tilde T}_J \cdots {\tilde T}_2{\tilde T}_1(c(\bp{p})).
\end{equation}
In this sense the canonical system essentially
grasps the time development pattern of the original one.
This fact, firstly recognized in \cite{TNS} for $\forall \theta_n = 1$,
will be utilized in Appendix \ref{app:a}.

\section{Combinatorial $R$ matrix as scattering matrix of
ultra-discrete solitons}\label{sec:proof}

Here we prove Theorem \ref{th:mainhomo}, which
identifies the scattering matrix of the
ultra-discrete solitons with the
combinatorial $R$ matrix of $U'_q(A^{(1)}_{M-1})$.

\subsection{Solitons}\label{subsec:solitons}

Let $B'_k$ be the classical crystal of $U'_q(A^{(1)}_{M-1})$
corresponding to the $k$-fold symmetric tensor representation:
\begin{equation*}
B'_k = \{\fbox{$m_1 \cd m_k$}  \mid m_i \in \{1,\ldots, M\},
m_1 \le \cd \le m_k \}.
\end{equation*}
Denote the Kashiwara operators of $B'_k$ by
$\ft'_i$ and $\et'_i$ for $i = 0,1,\ldots, M-1$.
For distinction, from now on we use the notation $B_k,\ft_i,\et_i$ for
$U'_q(A^{(1)}_M)$-crystals and
$B'_k,\ft'_i,\et'_i$ for $U'_q(A^{(1)}_{M-1})$-crystals.
Let $R$ and $R'$ be the combinatorial $R$ matrices for
$U'_q(A^{(1)}_M)$ and $U'_q(A^{(1)}_{M-1})$, respectively.
Thus $R \ft_i = \ft_i R$ and $R' \ft'_i = \ft'_i R'$ hold
when they act on the tensor product of two crystals,
and similarly for $\et_i, \et'_i$.
(We will specify the crystals that they act each time.)

\begin{remark}
When $M=1$ we still define $B'_k$ as above, which is
the set with the single element $u_k = \fbox{$1\ldots 1$}$.
We further understand  that the
``$U'_q(A^{(1)}_0)$" combinatorial $R$ matrix
$R': B'_l \ot B'_k \rightarrow B'_k \ot B'_l$ is given by
$R'(u_l\ot u_k) = u_k \ot u_l$.
\end{remark}
For each $k \in {\mathbb Z}_{\ge 1}$
define a map $\imath_k$ by
\begin{align*}
&\imath_k :\;\;\; B'_k  \qquad \qquad\;\longrightarrow
\quad\quad\quad\quad (B_1)^{\ot k}\\
&\quad\fbox{$m_1 \cd m_k$} \; \quad \quad \mapsto \; \;
 \fbox{$m_k+1$} \ot \cd \ot
\fbox{$m_1+1$}.
\end{align*}
Let $k_1, \ldots, k_N \in {\mathbb Z}_{\ge 1}$ and
$L_0, \ldots, L_N \in {\mathbb Z}_{\ge 0}$ for some $N \in {\mathbb Z}_{\ge 1}$.
In terms of $\imath_k$ we further introduce a map
\begin{equation*}
\imath^{(L_0,\ldots,L_N)}_{k_1, \ldots, k_N} :
B'_{k_1} \ot \cdots  \ot B'_{k_N} \longrightarrow
 (B_1)^{\ot L_0+ \cdots + L_N + k_1 + \cdots + k_N}
\end{equation*}
by
\begin{align*}
&\imath^{(L_0,\ldots,L_N)}_{k_1, \ldots, k_N}(b_1 \ot b_2 \ot \cdots \ot b_N) \\
&= \fbox{$1$}^{\ot L_0} \ot \imath_{k_1}(b_1) \ot \fbox{$1$}^{\ot L_1} \ot
\imath_{k_2}(b_2) \ot \cd \ot
\fbox{$1$}^{\ot L_{N-1}}\ot \imath_{k_N}(b_N) \ot \fbox{$1$}^{\ot L_N}.
\end{align*}
In particular $\imath_k = \imath^{(0,0)}_k$.
The map $\imath^{(L_0,\ldots,L_N)}_{k_1, \ldots, k_N}$ is injective.
For each $k \in {\mathbb Z}_{\ge 1}$  let $\varsigma_k$ denote the map
\begin{align*}
&\quad\varsigma_k :\;\;\; (B_1)^{\otimes k}
  \qquad \;\longrightarrow  \quad B_k\\
&\quad\fbox{$m_1$} \ot \cdots \ot \fbox{$m_k$} \; \mapsto \;
\fbox{$m'_1 \ldots m'_k$}
\end{align*}
where $1 \le m'_1 \le \cdots \le m'_k \le M+1$ are
just the re-ordering of $m_1, \ldots, m_k$ into the weakly increasing order.
We assume that ${\cal L} :=  \sum_n \theta_n$ is sufficiently large.
We set
\begin{equation}\label{eq:thetahat}
{\hat \theta} = (\cd \ot \varsigma_{\theta_{n}}
\ot \varsigma_{\theta_{n+1}}\ot \cd) :
B_1^{\ot {\cal L}} \rightarrow \;
\cd \ot B_{\theta_n} \ot B_{\theta_{n+1}}\ot \cd.
\end{equation}

For non-negative integers $L_0, \ldots, L_N$
such that ${\cal L} = L_0+ \cdots + L_N + k_1 + \cdots + k_N$,
denote by $\iota^{(L_0,\ldots,L_N)}_{k_1, \ldots, k_N}$
the  composition ${\hat \theta}\circ\imath^{(L_0,\ldots,L_N)}_{k_1, \ldots, k_N}$,
i.e.,
\begin{equation}\label{eq:Nsoliton}
\iota^{(L_0,\ldots,L_N)}_{k_1, \ldots, k_N}:
B'_{k_1} \ot \cdots \ot B'_{k_N}
\xrightarrow{\imath^{(L_0,\ldots,L_N)}_{k_1, \ldots, k_N}} B_1^{\ot {\cal L}}
\xrightarrow{{\hat \theta}} \quad
\cd \ot B_{\theta_n} \ot B_{\theta_{n+1}}\ot \cd.
\end{equation}
Suppose that the image is obtained from the element
$\cd \ot u_{\theta_n} \ot u_{\theta_{n+1}} \ot \cd$
by replacing only the isolated segments
$u_{\theta_{n_i}} \ot u_{\theta_{n_i+1}} \ot \cd \ot u_{\theta_{n'_i}}\, (n_i \le n'_i)$
with some
$b_{\theta_{n_i}} \ot \cd \ot b_{\theta_{n'_i}}
\in B_{\theta_{n_i}} \ot \cd \ot B_{\theta_{n'_i}}$  for $1 \le i \le N$.
Assume further that the interval is sufficiently large,
namely, $n_i - n'_{i-1} \gg \max(k_1, \ldots. k_N)$
for any $2 \le i \le N$.
In such a case we call the image of (\ref{eq:Nsoliton})
an {\em asymptotic $N$ soliton state}.
Each soliton is essentially associated with an element in $B'_k$,  and
we call $k$ the {\em amplitude} of the corresponding soliton.
States obtained from an asymptotic $N$ soliton state under
arbitrary time evolutions $T_{\kappa} \cdots T_{\kappa'}$
will be called  $N$ soliton states.
This definition will naturally be justified from the consideration on the
conserved quantities in Section \ref{subsec:conserve}.
Note  that
$\iota^{(L_0,\ldots,L_N)}_{k_1, \ldots, k_N}$ is not injective
since ${\hat \theta}$ is not.
Consequently, the result of application of
$\iota^{(L'_0,\ldots,L'_N)}_{k_1, \ldots, k_N}$
is not necessarily an `overall translation' of (\ref{eq:Nsoliton}) in a naive sense
even when $L'_i - L_i$ is $i$-independent for $i < N$ or $i > 0$.
See Example \ref{ex:onesoltoki} below.

First we consider $N=1$ case.
As it turns out in Proposition \ref{pr:onesoliton},
there is no distinction between an
asymptotic 1 soliton state and a 1 soliton state.
Moreover one can check that the definition of the 1 soliton state here
agrees with the 1 soliton solution that will be given later in (\ref{1-sol.sol}).
Given a 1 soliton state
\begin{equation*}
\bp{p} = \cd \ot b_{n-1} \ot b_n \ot b_{n+1} \ot \cd \in
\cd \ot B_{\theta_{n-1}} \ot B_{\theta_n} \ot B_{\theta_{n+1}} \ot \cd,
\end{equation*}
one can unambiguously specify integers $n, k(\ge 1), s, t$
by the conditions:
\begin{align*}
&b_j = u_{\theta_j} \quad \text{ if } j < n \text{ or } j > n+k,\\
&b_n = \fbox{$1\cd 1 m_1 \cd m_t$} \quad 1 \le t \le \theta_n,
\, 2 \le m_1 \le \cd \le m_t \le M+1,\\
&b_{n+k} = \fbox{$1\cd 1 m'_1 \cd m'_s$} \quad 0 \le s \le \theta_{n+k}-1,
\, 2 \le m'_1 \le \cd \le m'_s \le M+1,\\
&b_{n+1}, \ldots, b_{n+k-1} \text{ do not contain } 1 \text{ in their tableaux}.
\end{align*}
Note that `if' in the first condition is not `only if' in that
$b_{n+k}=u_{\theta_{n+k}}$ is allowed as $s=0$.
The amplitude of the soliton according to the above definition equals
$t + \theta_{n+1} + \cd + \theta_{n+k-1} + s$.
We set
\begin{equation*}
x(\bp{p}) = \sum_{j\le n} \theta_j - t, \qquad
y(\bp{p}) = t + \theta_{n+1} + \cd + \theta_{n+k-1}
\end{equation*}
and call $x(\bp{p})$ the {\em coordinate} of the soliton.
$y(\bp{p})$ should not be confused with the amplitude of the soliton.
\begin{example}\label{ex:onesoltoki}
Consider $B_{\theta_1} \ot \cd \ot B_{\theta_6}$ with
$\theta_1=\theta_3=1, \, \theta_2=\theta_4=\theta_6=2$ and $\theta_5=3$, hence
${\cal L} = 11$.
\par
(i) Take $b = \fbox{$1$} \in B'_1$.
Then $\iota^{(L_0,L_1)}_1(b)$ with
$L_0 + L_1 = 10$ are
examples of 1 soliton states with amplitude 1.
One has $\iota^{(1,9)}_1(b) =
\iota^{(2,8)}_1(b)$,
$\iota^{(4,6)}_1(b) =
\iota^{(5,5)}_1(b)$ and
$\iota^{(6,4)}_1(b) =
\iota^{(7,3)}_1(b) =
\iota^{(8,2)}_1(b)$.
For $L_0 \le 8$ they look as
\small{
\begin{align*}
\bp{p}\quad&\qquad\qquad\qquad\qquad\quad\qquad \qquad\qquad\qquad\;
n\quad n+k\quad x(\bp{p})\quad y(\bp{p})\\
\iota^{(0,10)}_1(b)
&= \fbox{$2$}\ot\fbox{$11$}\ot\fbox{$1$}\ot\fbox{$11$}
\ot\fbox{$111$}\ot\fbox{$11$}
\quad  1\quad\quad 2\, \quad \;\quad 0 \quad\;\quad 1,\\
\iota^{(2,8)}_1(b)
 &= \fbox{$1$}\ot\fbox{$12$}\ot\fbox{$1$}\ot\fbox{$11$}
\ot\fbox{$111$}\ot\fbox{$11$}
\quad  2\quad\quad 3\, \quad \;\quad 2 \quad\;\quad 1,\\
\iota^{(3,7)}_1(b)
&= \fbox{$1$}\ot\fbox{$11$}\ot\fbox{$2$}\ot\fbox{$11$}
\ot\fbox{$111$}\ot\fbox{$11$}
\quad  3\quad\quad 4\, \quad \;\quad 3 \quad\;\quad 1,\\
\iota^{(5,5)}_1(b)
 &= \fbox{$1$}\ot\fbox{$11$}\ot\fbox{$1$}\ot\fbox{$12$}
\ot\fbox{$111$}\ot\fbox{$11$}
\quad  4\quad\quad 5\, \quad \;\quad 5 \quad\;\quad 1,\\
\iota^{(8,2)}_1(b)
 &= \fbox{$1$}\ot\fbox{$11$}\ot\fbox{$1$}\ot\fbox{$11$}
\ot\fbox{$112$}\ot\fbox{$11$}
\quad  5\quad\quad 6\, \quad \;\quad 8 \quad\;\quad 1,
\end{align*}
}
where we have also listed $n, n+k, x(\bp{p})$ and $y(\bp{p})$.
\par
(ii) Take $c = \fbox{$11223$} \in B'_5$.
Then
$\iota^{(L_0,L_1)}_1(b)$ with
$L_0 + L_1 = 6$ are
examples of 1 soliton states with amplitude 5.
For $L_0 \le 5$ they look as
\small{
\begin{align*}
\bp{p}\quad&\qquad\qquad\qquad\qquad\quad\qquad \qquad\qquad\qquad\;
n\quad n+k\quad x(\bp{p})\quad y(\bp{p})\\
\iota^{(0,6)}_5(c) &= \fbox{$4$}\ot\fbox{$33$}\ot\fbox{$2$}\ot\fbox{$12$}
\ot\fbox{$111$}\ot\fbox{$11$}
\quad  1\quad\quad 4\, \quad \;\quad 0 \quad\;\quad 4,\\
\iota^{(1,5)}_5(c) &= \fbox{$1$}\ot\fbox{$34$}\ot\fbox{$3$}\ot\fbox{$22$}
\ot\fbox{$111$}\ot\fbox{$11$}
\quad  2\quad\quad 5\, \quad \;\quad 1 \quad\;\quad 5,\\
\iota^{(2,4)}_5(c) &= \fbox{$1$}\ot\fbox{$14$}\ot\fbox{$3$}\ot\fbox{$23$}
\ot\fbox{$112$}\ot\fbox{$11$}
\quad  2\quad\quad 5\, \quad \;\quad 2 \quad\;\quad 4,\\
\iota^{(3,3)}_5(c) &= \fbox{$1$}\ot\fbox{$11$}\ot\fbox{$4$}\ot\fbox{$33$}
\ot\fbox{$122$}\ot\fbox{$11$}
\quad  3\quad\quad 5\, \quad \;\quad 3 \quad\;\quad 3,\\
\iota^{(4,2)}_5(c) &= \fbox{$1$}\ot\fbox{$11$}\ot\fbox{$1$}\ot\fbox{$34$}
\ot\fbox{$223$}\ot\fbox{$11$}
\quad  4\quad\quad 6\, \quad \;\quad 4 \quad\;\quad 5,\\
\iota^{(5,1)}_5(c) &= \fbox{$1$}\ot\fbox{$11$}\ot\fbox{$1$}\ot\fbox{$14$}
\ot\fbox{$233$}\ot\fbox{$12$}
\quad  4\quad\quad 6\, \quad \;\quad 5 \quad\;\quad 4.
\end{align*}
}
\end{example}
In Section \ref{subsec:2solgen} we will make use of
\begin{proposition}\label{pr:onesoliton}
Let $\bp{p} = \iota^{(L_0, L_1)}_l(b)$ be the
1 soliton of amplitude $l$ associated with $b \in B'_l$.
Then its time evolution  $T_\kappa(\bp{p})$  is
again 1 soliton and  expressible as
$T_\kappa(\bp{p}) = \iota^{(L'_0, L'_1)}_l(b)$ for some
$L'_0, L'_1$ $ (L'_0 + L'_1 = L_0 + L_1)$ but with the same $b \in B'_l$.
The difference of their coordinates (velocity under $T_\kappa$) is given by
\begin{equation*}
x(T_\kappa(\bp{p})) - x(\bp{p}) =
\begin{cases}
\kappa & \kappa < y(\bp{p}), \\
\min(\kappa,l) + \max(\theta_{n+k}-l, 0) & \kappa \ge y(\bp{p}).
\end{cases}
\end{equation*}
\end{proposition}
The proof is done by a cumbersome case study.
When $\forall \theta_n = 1$, the above result simplifies to
$x(T_\kappa(\bp{p})) - x(\bp{p}) = \min(\kappa, l)$
in agreement with \cite{FOY}.
In general, the velocity varies locally depending on the data $\{ \theta_n \}$.
In Example \ref{ex:onesoltoki} (i) one has
$T_\kappa(\iota^{(0,10)}_1(b))
= \iota^{(2,8)}_1(b),
T_\kappa(\iota^{(2,8)}_1(b))
= \iota^{(3,7)}_1(b),
T_\kappa(\iota^{(3,7)}_1(b))
= \iota^{(5,5)}_1(b),
T_\kappa(\iota^{(5,5)}_1(b))
= \iota^{(8,2)}_1(b)$ for any $\kappa \ge 1$.
Similarly in (ii) one has
$T_\kappa(\iota^{(0,6)}_5(c))
 = \iota^{(\kappa',6-\kappa')}_5(c)$
for any $\kappa \ge 1$, where $\kappa' = \min(\kappa,5)$.
These results agree with Proposition \ref{pr:onesoliton}.

Let $\iota^{(L_0,\ldots, L_N)}_{k_1, \ldots, k_N}(c_1 \ot \cd \ot c_N)
\, (c_i \in B'_{k_i})$
be an asymptotic $N$ soliton state and
\begin{equation*}
\cd \ot b^t_n \ot b^t_{n+1} \ot \cd =
T_{\kappa_t}T_{\kappa_{t-1}} \cdots
\left(\iota^{(L_0,\ldots, L_N)}_{k_1, \ldots, k_N}
(c_1 \ot \cd \ot c_N)
\right)
\end{equation*}
be its time evolution.
Assume that the solitons are enough separated without an interaction
throughout the time interval in consideration.
Let $\{ v^t_n\}$ be the associated variables on the vertical edges
as in (\ref{eq:vbbv}).
Then in the space-time interchanged picture, the state
$\cd \ot v^{t+1}_n \ot v^t_n \ot \cd$ is also an asymptotic $N$ soliton state
associated with the same $c_1 \ot \cd \ot c_N$.
Namely,
\begin{equation*}
\cd \ot v^{t+1}_n \ot v^t_n \ot \cd =
{\hat \kappa} \circ \imath^{(L'_0,\ldots, L'_N)}_{k_1, \ldots, k_N}
(c_1 \ot \cd \ot c_N)
\end{equation*}
for some $L'_0,\ldots, L'_N$.
Here
\begin{equation*}
{\hat \kappa} = (\cd \ot \varsigma_{\kappa_{t+1}}
\ot \varsigma_{\kappa_{t}}\ot \cd) :
B_1^{\ot {\cal M}} \rightarrow \;
\cd \ot B_{\kappa_{t+1}} \ot B_{\kappa_{t}}\ot \cd.
\end{equation*}
is an analogue of ${\hat \theta}$ in (\ref{eq:thetahat}), and
we have set ${\cal M} = \sum_t \kappa_t$.
The figure in Example \ref{ex:tnsd} will be of help to understand this fact.
In a sense one can employ either picture to
describe the scattering process.
Indeed our discussion in the end of Section \ref{subsec:2solgen}
will rely on  this observation.


\subsection{Scattering of 2 solitons; a typical case}
\label{subsec:2solitons}

Our aim here is to show Theorem \ref{th:main48} which is valid in
the `typical' situation (\ref{eq:assumptoki}).
\begin{lemma}\label{lem:lem1}
For each $i = 1, \ldots, M-1$, we have a commutative diagram:
\begin{equation*}
\begin{CD}
B'_k @>{\imath_k}>> (B_1)^{\ot k}\\
@V{\tilde{e}'_i}VV @VV{\et_{i+1}}V\\
B'_k \sqcup\{0\}@>{\imath_k}>> (B_1)^{\ot k}\sqcup\{0\},
\end{CD}
\end{equation*}
where $\imath_k(0) = 0$.
The same relation holds also between  $\ft'_i$ and $\ft_{i+1}$.
\end{lemma}
Combining Lemma \ref{lem:lem1} with the
realization of $B_\theta$ in $B_1^{\ot \theta}$
as a $U_q(A_M)$-crystal (cf. \cite{KN}),
one can derive the following lemmas.
\begin{lemma}\label{lem:cor1}
In the diagram
\begin{equation*}
\begin{CD}
B'_l \ot B'_k @>{\iota^{(L_0,L_1,L_2)}_{l,k}}>>
(\cd \ot B_{\theta_n} \ot B_{\theta_{n+1}}\ot  \cd) \\
@V{\tilde{e}'_i}VV @VV{\et_{i+1}}V\\
(B'_l \ot B'_k) \sqcup\{0\}@>{\iota^{(L_0,L_1,L_2)}_{l,k}}>>
(\cd \ot B_{\theta_n} \ot B_{\theta_{n+1}}\ot  \cd)\sqcup\{0\},
\end{CD}
\end{equation*}
suppose that the image of $\iota^{(L_0,L_1,L_2)}_{l,k}$
is an asymptotic 2 soliton state.
Then the diagram is commutative for any $i = 1, \ldots, M-1$.
The same relation holds also between  $\ft'_i$ and $\ft_{i+1}$.
\end{lemma}
Actually, the commutativity of the above diagram holds
under a milder condition than being an asymptotic  2 soliton state.
\begin{lemma}\label{lem:lem2}
Let $p_1,\ldots, p_m$
be the subsequence of
$a_1, \ldots, a_L\,  (a_n \in B_{\theta_n})$ consisting of all the elements
such that $a_n \neq u_{\theta_n}$.
Assume the same relation between
$p'_1, \ldots, p'_m$ and
$a'_1, \ldots, a'_L$.
Then for any $t, t' \in {\mathbb Z}_{\ge 0}$ and
$\kappa \in {\mathbb Z}_{\ge 1}$, the two relations
\begin{align*}
&\ft_{i+1}(p_1 \ot \cd \ot p_m ) = p'_1 \ot \cd \ot p'_m,\\
&\ft_{i+1}(u_\kappa^{\ot t} \ot a_1 \ot \cd \ot a_L \ot u_\kappa^{\ot t'}) =
u_\kappa^{\ot t} \ot a'_1 \ot \cd \ot a'_L \ot u_\kappa^{\ot t'}
\end{align*}
are equivalent
for each $i = 1, \ldots, M-1$.
The equivalence persists even when the right hand sides are both $0$.
The same is true also for $\et_{i+1}$.
\end{lemma}
\begin{proposition}\label{prop:sf=fs}
Suppose an asymptotic two soliton state has evolved into another as
\begin{equation}\label{eq:scatter}
T^t_\kappa \left( \iota^{(L_0,L_1,L_2)}_{l,k}(b \ot c) \right)
= \iota^{(L'_0,L'_1,L'_2)}_{k,l}(c' \ot b')
\end{equation}
for some $\kappa, t, L_i, L'_i > 0, b, b' \in B'_l$ and $c, c' \in B'_k$.
Then (\ref{eq:scatter}) is also valid under the replacement of
$b \ot c$ (resp. $c' \ot b'$) by
$\ft'_i(b \ot c)$ (resp. $\ft'_i(c' \ot b')$) for any $i = 1, \ldots, M-1$
such that $\ft'_i(b \ot c) \neq 0$.
\end{proposition}
{\em Proof}.
(\ref{eq:scatter}) is equivalent to
\begin{equation*}
u_\kappa^{\ot t}\ot
\iota^{(L_0,L_1,L_2)}_{l,k}(b \ot c)
\simeq  \iota^{(L'_0,L'_1,L'_2)}_{k,l}(c' \ot b')
\ot u_\kappa^{\ot t}
\end{equation*}
Apply $\ft_{i+1}$ to the both sides.
Due to Lemmas \ref{lem:cor1} and \ref{lem:lem2}, the result becomes
\begin{equation*}
u_\kappa^{\ot t}\ot
\iota^{(L_0,L_1,L_2)}_{l,k}(\ft'_i(b \ot c))
\simeq  \iota^{(L'_0,L'_1,L'_2)}_{k,l}(\ft'_i(c' \ot b'))
\ot u_\kappa^{\ot t}.
\end{equation*}
\begin{flushright}
$\square$
\end{flushright}
\par

\begin{proposition}\label{prop:toki}
Let   $l > k$ and assume that
$\iota^{(L_0,L_1,L_2)}_{l,k}(b_1 \ot b_2)$
is an asymptotic  2 soliton state with
\begin{equation}\label{eq:b1b2}
b_1 = (l,0, \ldots, 0) \in B'_l, \quad
b_2 = (h,k-h,0, \ldots, 0) \in B'_k
\end{equation}
with $0 \le h \le k$ in the notation of (\ref{eq:action}).
Assume further that $l > \theta_n$ for all but finitely many $n$'s.
Then if $\kappa \gg l$, there exists $t > 0$ such that the result of the time evolution
$T^t_\kappa$ also becomes the asymptotic 2 soliton state as
\begin{equation}\label{eq:hscatter}
T^t_\kappa \left( \iota^{(L_0,L_1,L_2)}_{l,k}(b_1 \ot b_2) \right)
= \iota^{(L'_0,L'_1,L'_2)}_{k,l}(c_2 \ot c_1),
\end{equation}
where $c_1, c_2$ are given by
\begin{equation}\label{eq:c2c1}
c_2 = (k,0, \ldots,0) \in B'_k, \quad
c_1 = (l-k+h,k-h,0, \ldots,0) \in B'_l.
\end{equation}
\end{proposition}
The proof is given in Appendix \ref{app:a}.
In fact both $b_1 \ot b_2$ and $c_2 \ot c_1$ are $U_q(A_{M-1})$
highest element, i.e.,
$\et'_i(b_1\ot b_2) = \et'_i(c_2\ot c_1) = 0$ for all $1 \le i \le M-1$.
Combining this property with the
conservation of  weights (number of the letters)
and the soliton content (cf. Section \ref{subsec:conserve}),
one can argue that the
outgoing state should necessarily correspond to
$c_2 \ot c_1$ if it is an asymptotic 2 soliton state at all.
However, to establish the separation into 2 solitons asymptotically
is not a trivial task for  inhomogeneous $\theta_n$'s
only bounded by the condition $l \ge \theta_n$
for all but finitely many $n$'s.
So far we have not managed it without recourse to
the actual 2 soliton solution as in Appendix \ref{app:a}.
\par

As a $U_q(A_{M-1})$-crystal, the $U_q'(A^{(1)}_{M-1})$-crystal
$B'_l \ot B'_k$ decomposes into the connected components.
Each component is parametrized with the $U_q(A_{M-1})$
highest elements $b_1 \ot b_2$ (\ref{eq:b1b2}), and is
generated by applying $\ft'_i$ operators ($1 \le i \le M-1$) to it.
The decomposition of the same pattern takes place also for $B'_k \ot B'_l$
according to the highest elements $c_2 \ot c_1$.
Combining this fact with Propositions \ref{prop:sf=fs} and
\ref{prop:toki}, we conclude that
there exists a map $S'$ ($S$ matrix) uniquely defined by
\begin{align}\label{eq:sdef}
&S': \;\; B'_l \otimes B'_k \,\rightarrow \,B'_k \otimes B'_l\nonumber\\
&T^t_\kappa \left( \iota^{(L_0,L_1,L_2)}_{l,k}(b \ot c) \right)
= \iota^{(L'_0,L'_1,L'_2)}_{k,l}(S'(b \ot c)),
\end{align}
under the condition
\begin{equation}\label{eq:assumptoki}
\kappa \gg l > k, \quad  l > \theta_n \, \text{ for all but finitely many } \,
n\text{'s}.
\end{equation}
It describes the 2 soliton scattering.

\begin{theorem}\label{th:main48}
Under the assumption (\ref{eq:assumptoki}), we have
$R' = S'$ on the $U'_q(A^{(1)}_{M-1})$-crystal
$B'_l \ot B'_k$.
\end{theorem}
{\em Proof}.
By the definition and Proposition \ref{prop:sf=fs},
the both $R'$ and $S'$ commute with $\ft'_i$ for any
$1 \le i \le M-1$.
Moreover, for any $U_q(A_{M-1})$ highest elements $b_1 \ot b_2$ given by
(\ref{eq:b1b2}),
their actions are the same, i.e.,
$S'(b_1 \ot b_2) = c_2 \ot c_1 = R'(b_1 \ot b_2)$,
where the latter can be verified {}from Proposition \ref{prop:H-rule}.
\hfill$\square$

Thus in the situation (\ref{eq:assumptoki}) the larger soliton
overtakes the smaller soliton and the scattering matrix coincides with
the combinatorial $R$ matrix of $U'_q(A^{(1)}_{M-1})$-crystal.
For instance Example \ref{ex:tns} tells that
\[
S': \fbox{$13$} \ot \fbox{$2$} \mapsto \fbox{$1$} \ot \fbox{$23$} \,.
\]
This agrees with Example \ref{ex:R} (i).

\subsection{Scattering of 2 solitons; general case}\label{subsec:2solgen}

First let us consider the homogeneous case
$\forall \theta_n = \theta, \forall \kappa_t = \kappa$.
Fix positive integers $l > k$.
We study the scattering of 2 solitons in
$\cd \ot B_\theta \ot B_\theta \ot \cd$
with amplitudes $l$ and $k$ under the time evolution
$T^t_\kappa$.
The qualitative feature of the scattering depends on the cases:
\begin{alignat*}{2}
&(\text{i})\; l > k \ge \max(\theta, \kappa) & \quad v_l = v_k = \kappa,\\
&(\text{ii}) \;  \min(\theta, \kappa) \ge l > k & \quad v_l = v_k = \theta,\\
&(\text{iii}) \; l \ge \kappa > k \ge \theta & \quad v_l = \kappa > v_k = k,\\
&(\text{iv}) \; \kappa \ge l > k \ge \theta & \quad v_l = l > v_k = k,\\
&(\text{v}) \; l \ge \kappa > \theta \ge k & \quad v_l = \kappa > v_k = \theta,\\
&(\text{vi}) \; \kappa \ge l > \theta \ge k & \quad v_l = l > v_k = \theta,\\
&(\text{vii}) \; l \ge \theta > k \ge \kappa, & \\
&(\text{viii}) \; \theta \ge l > k \ge \kappa, & \\
&(\text{ix}) \; l \ge \theta > \kappa \ge k, & \\
&(\text{x}) \; \theta \ge l > \kappa \ge k. &
\end{alignat*}
Here the classification has been done so that
\begin{equation*}
\left\{ (\text{i}) \amalg (\text{ii}) \right\} \coprod
\left\{ (\text{iii}) \cup (\text{iv}) \cup (\text{v}) \cup (\text{vi}) \right\} \coprod
\left\{ (\text{vii}) \cup (\text{viii}) \cup (\text{ix}) \cup (\text{x}) \right\}.
\end{equation*}
For example (iii) and (iv) share $l=\kappa > k \ge \theta$ case.
However the three groups are mutually disjoint and correspond to
distinct features of the scattering as we will see below.
The $v_l$ and $v_k$ are the velocities of the solitons with amplitude $l$ and $k$,
respectively.
For each soliton it has been calculated by using Proposition \ref{pr:onesoliton}
by assuming no effect from the other soliton.
In (vi) and (x)
we have excluded  $l = \theta$ and $l = \kappa$, respectively
since they both lead to  $v_l = v_k = \theta$ hence no scattering.
By the same reason the cases  (i) and (ii) are out of question.
Via the space-time interchange $\theta \leftrightarrow \kappa$,
the cases (vii),(viii), (ix) and (x) are mapped to
(iii), (iv), (v) and (vi), respectively.
(See the argument before Theorem \ref{th:mainhomo} on the
velocities in the cases (vii) - (x).)
Thus we are left with (iii)-(vi), where
$l > \theta$ and $v_l > v_k$ are always valid.
Following \cite{FOY}, we utilize the commutativity (\ref{eq:commute})
and  consider the 2 soliton scattering under $T^t_\kappa$ as
\begin{equation*}
T^t_\kappa = T^{-t'}_\infty \,T^t_\kappa \,T^{t'}_\infty.
\end{equation*}
The scattering are thus divided into three stages.
In the first stage, we let solitons evolve under $T^{t'}_\infty$
for sufficiently large $t'$.
Since $l > \theta$ matches the condition (\ref{eq:assumptoki}),
Theorem \ref{th:main48} tells that the larger soliton overtakes
the smaller one with the scattering rule described by $S' = R'$.
In the second stage corresponding to $T^t_\kappa$,
the larger soliton goes further ahead than the smaller one
with no interaction because of $v_l > v_k$.
Therefore in the last stage $T^{-t'}_\infty$,
the two remain isolated even though they are drawn back and
get relatively closer.
Thus we conclude that in all the cases (iii)-(vi),
the qualitative feature is the same as the one in Theorem \ref{th:main48}.
Namely, the larger soliton overtakes the smaller one and the
scattering rule is given by the combinatorial $R$ matrix
$R': B'_l \ot B'_k \rightarrow B'_k \ot B'_l$.
Through the space-time interchange argument,
this implies the opposite feature of scattering in the cases (vii)-(x).
Namely, the smaller one overtakes the larger one with the scattering rule
given by the combinatorial $R$ matrix
$R': B'_k \ot B'_l \rightarrow B'_l \ot B'_k$.

We note that in  the cases (vii) - (x),
one does not necessarily have $v_l < v_k$ at any time.
It actually
depends on whether $\kappa \gtrless y(\bp{p})$
when Proposition \ref{pr:onesoliton} is applied.
Nevertheless $v_l < v_k$ should be valid ``on average" and the
above feature of the scattering should hold  due to
the reduction to the cases (iii) - (vi) where
the strict inequality $v_l > v_k$ is always valid.
To summarize, we have shown
\begin{theorem}\label{th:mainhomo}
Let $l > k$ be the amplitude of 2 solitons in
$\cd \ot B_\theta \ot B_\theta \ot \cd$.
Under the time evolution $T_\kappa$,
the scattering matrix of the collision (if any) in the sense of
(\ref{eq:sdef}) or (\ref{eq:sdef}$)_{l\leftrightarrow k}$ is
given by $S' = R'$, where $R'$ is the combinatorial $R$
matrix of the $U'_q(A^{(1)}_{M-1})$-crystals for
\begin{align*}
(\text{I})& \; B'_l \ot B'_k \stackrel{\sim}{\mapsto} B'_k \ot B'_l
 \; \text{ if } \min(l,\kappa) > \max(k,\theta),\\
(II)& \; B'_k \ot B'_l \stackrel{\sim}{\mapsto} B'_l \ot B'_k
 \; \text{ if } \min(l,\theta) > \max(k,\kappa),\\
(III)& \; \text{no scattering (same velocity)} \; \text{ otherwise}.
\end{align*}
\end{theorem}

Example \ref{ex:tnsd} corresponds to the choice
$l = \theta = 2, k = \kappa = 1$, hence to (II) in the theorem.
The scattering matrix is read off the figure,
\[
S': \fbox{$3$} \ot \fbox{$22$} \mapsto \fbox{$23$} \ot \fbox{$2$} \,.
\]
This agrees with the inverse of the $R$ matrix in  Example \ref{ex:R} (ii).

Let us comment on the inhomogeneous case where
$\theta_n$'s and $\kappa_t$'s actually depend on the indices.
In view of (\ref{eq:assumptoki}),
the qualitative feature of the scattering remains the same
as Theorem \ref{th:mainhomo} even if we slightly relax the conditions therein.
For instance  the larger soliton still overtakes the smaller one with the rule
$S' = R'$ if $\min(l,\kappa_t) > \max(k,\theta_n)$ holds for almost all $n$ and $t$
that are relevant during the scattering in question.
In such cases we expect that
the asymptotic $N$ soliton state in the sense of Section \ref{subsec:solitons}
undergoes the scattering which are essentially factorized into the two-body ones
studied here.
On the other hand, if $\theta_n$'s and $\kappa_t$'s are not
bounded by the condition as above and indeed far from being homogeneous,
even 2 solitons can collide many times in general depending on the local velocities.
In such a case we do not have a simple picture of the scattering.

\begin{example}\label{ex:doublescattering}
Let $M=3$.
\small{
\begin{center}
$0:\cdots 14 \cdot 3 \cdot 123 \cdot 111 \cdot 24
   \cdot 1 \cdot 1 \cdot 111 \cdot 11  \cdot 1
\cdot 1111 \cdot 1111 \cdot 11111 \cdot 11111 \cdot 111 \cdot 1111\cdots$\\
$1:\cdots 11 \cdot 1 \cdot 114 \cdot 233 \cdot 11
   \cdot 4 \cdot 2 \cdot 111 \cdot 11  \cdot 1
\cdot 1111 \cdot 1111 \cdot 11111 \cdot 11111 \cdot 111 \cdot 1111\cdots$\\
%
$2:\cdots 11 \cdot 1 \cdot 111 \cdot 111 \cdot 34
   \cdot 3 \cdot 1 \cdot 224 \cdot 11  \cdot 1
\cdot 1111 \cdot 1111 \cdot 11111 \cdot 11111 \cdot 111 \cdot 1111\cdots$\\
%
$3:\cdots 11 \cdot 1 \cdot 111 \cdot 111 \cdot 11
   \cdot 1 \cdot 4 \cdot 113 \cdot 34  \cdot 2
\cdot 1112 \cdot 1111 \cdot 11111 \cdot 11111 \cdot 111 \cdot 1111\cdots$\\
%
$4:\cdots 11 \cdot 1 \cdot 111 \cdot 111 \cdot 11
   \cdot 1 \cdot 1 \cdot 114 \cdot 13  \cdot 1
\cdot 1234 \cdot 1112 \cdot 11111 \cdot 11111 \cdot 111 \cdot 1111\cdots$\\
%
$5:\cdots 11 \cdot 1 \cdot 111 \cdot 111 \cdot 11
   \cdot 1 \cdot 1 \cdot 111 \cdot 14  \cdot 3
\cdot 1114 \cdot 1223 \cdot 11111 \cdot 11111 \cdot 111 \cdot 1111\cdots$\\
%
$6:\cdots 11 \cdot 1 \cdot 111 \cdot 111 \cdot 11
   \cdot 1 \cdot 1 \cdot 111 \cdot 11  \cdot 1
\cdot 1134 \cdot 1234 \cdot 11112 \cdot 11111 \cdot 111 \cdot 1111\cdots$\\
%
$7:\cdots 11 \cdot 1 \cdot 111 \cdot 111 \cdot 11
   \cdot 1 \cdot 1 \cdot 111 \cdot 11  \cdot 1
\cdot 1111 \cdot 2334 \cdot 11114 \cdot 11112 \cdot 111 \cdot 1111\cdots$\\
%
$8:\cdots 11 \cdot 1 \cdot 111 \cdot 111 \cdot 11
    \cdot 1 \cdot 1 \cdot 111 \cdot 11  \cdot 1
\cdot 1111 \cdot 1134 \cdot 11123 \cdot 11114 \cdot 112 \cdot 1111\cdots$\\
%
$9:\cdots 11 \cdot 1 \cdot 111 \cdot  111  \cdot  11
    \cdot 1 \cdot 1 \cdot 111 \cdot 11  \cdot 1
\cdot 1111  \cdot  1111  \cdot 12334  \cdot 11111  \cdot 114  \cdot 1112\cdots$
\end{center}
}\noindent
where $\cdot$ denotes $\otimes$, and $14$
for example does \fbox{$14$} $\in B_2$.
Not only $\theta_n$'s but also $\kappa_t$ are inhomogeneous here
so that the relevant time evolutions are
$T_5$ for the process
$0 \rightarrow 1 \rightarrow 2 \rightarrow 3 \rightarrow 4$, whereas they are
$T_2$ for
$4 \rightarrow 5 \rightarrow 6 \rightarrow 7 \rightarrow 8 \rightarrow 9$.
This is an example of the double scattering of 2 solitons
caused by the inhomogeneity.
The larger soliton once overtakes the smaller one, but
after the collision it gets slower due to the
environmental change and is eventually passed by the smaller one again.
This is easily understood from
the classification (I)--(III) in Theorem \ref{th:mainhomo}
for the homogeneous case.
In the first stage we have $l=4, \kappa=5,
k=2, \theta_n \le 3$ so that the larger soliton overtakes the smaller as in (I).
On the other hand we have
$\kappa=2, \theta_n \ge 3$ in the second stage
hence the smaller one passes the larger one as  in (II).
Following the time evolution downward,
one finds the scattering matrices for the successive collisions:
\begin{equation*}
\fbox{$1223$} \otimes \fbox{$13$} \mapsto
\fbox{$23$} \otimes \fbox{$1123$} \mapsto
\fbox{$1223$} \otimes \fbox{$13$}
\end{equation*}
in terms of the soliton labels with the $U'_q(A^{(1)}_2)$-crystal elements.
They agree with the combinatorial $R$ matrices
$B'_4 \ot B'_2 \simeq B'_2 \ot B'_4$
calculated from Proposition \ref{prop:H-rule}.
\end{example}

\subsection{Conserved quantities}\label{subsec:conserve}

Let us give a class of conserved quantities in the $A^{(1)}_M$
automaton.
Since our construction here is based on \cite{FOY} and the result is
quite parallel, we will only present a brief sketch.
Given an automaton state
$\mbox{\mathversion{bold}$p$} = \cdots \ot b_{n} \ot
b_{n+1} \ot \cdots \, (b_n = u_{\theta_n} \text{ for } \vert n \vert \gg 1)$,
let
\begin{equation*}
u_{\kappa} \ot \mbox{\mathversion{bold}$p$} \simeq
\cd \ot b'_{n-1}\ot b'_n \ot v_n \ot b_{n+1} \ot b_{n+2} \ot \cd
\end{equation*}
for some $b'_i \in B_{\theta_i}$ and $v_n \in B_\kappa$.
Set
\begin{equation*}
E_\kappa(\mbox{\mathversion{bold}$p$}) =
-\sum_n H_{B_\kappa B_{\theta_{n+1}}}(v_n \ot b_{n+1}),
\end{equation*}
which is well defined owing to the normalization (\ref{eq:Hnormalize}).
By the same argument as in \cite{FOY} we get
\begin{equation*}
E_\kappa(T_{\kappa'}(\mbox{\mathversion{bold}$p$})) =
E_\kappa(\mbox{\mathversion{bold}$p$}) \quad \text{ for any } \kappa, \kappa'.
\end{equation*}
Thus $E_\kappa(\mbox{\mathversion{bold}$p$}), \, \kappa \in {\mathbb Z}_{\ge 1}$
form a family of conserved quantities.
If $\mbox{\mathversion{bold}$p$}$ is an asymptotic
$N$ soliton state in the sense of Section \ref{subsec:solitons},
it is straightforward to derive
\begin{equation}\label{eq:evalue}
E_\kappa(\mbox{\mathversion{bold}$p$}) = \sum_{l\ge 1} \min(l,\kappa)N_l,
\end{equation}
where $N_l$ is the number of solitons with amplitude $l$.
Therefore if a state with the soliton content
$\{N_l\}$ scatter into another state with the content
$\{N'_l\}$, $N_l = N'_l$ must be valid for any $l$
due to the conservation of all $E_\kappa$'s.
In both Example \ref{ex:tns} and \ref{ex:tnsd}, we have
$E_1 = 2, E_l= 3$ for $l \ge 3$, in agreement  with $N_1=N_2=1, N_l = 0$
for $l \ge 3$.
In Example \ref{ex:doublescattering}, we have
$E_1=2, E_2=4, E_3=5$ and $E_l=6$ for $l \ge 4$,
in agreement with $N_l = \delta_{l 2}+\delta_{l 4}$.

When  $\forall \theta_n=1$, (\ref{eq:evalue}) is obtained
in proposition 4.4 in \cite{FOY}.
An equivalent family of the conserved quantities has
also been given in \cite{TNS}.

Another conserved quantity is the semistandard Young tableau, which can be
constructed as follows.
Given an automaton state
$\mbox{\mathversion{bold}$p$} = \cdots \ot b_{n} \ot
b_{n+1} \ot \cdots$, let
$c_s \ldots c_2 c_1$ be the subsequence of
$\ldots b_{n-1} b_n b_{n+1} \ldots$ obtained by dropping all the
$b_j$'s such that $b_j = u_{\theta_j}$.
Each $c_j$ has the form
\begin{equation*}
c_j = \fbox{$1\ldots 1 m_1 \ldots m_k$} \quad 2 \le m_1 \le \cd \le m_k \le M+1\;
\text{ for some } k \ge 1,
\end{equation*}
for which we set
\begin{equation*}
\overline{c}_j = \fbox{$m_1-1 \ldots m_k-1$} \in B'_k.
\end{equation*}
Let ${\cal T}(\mbox{\mathversion{bold}$p$})
:= (((\overline{c}_1  \cdot \overline{c}_2) \cdot
\overline{c}_3) \cdot \; \cdots \; \cdot \overline{c}_s)$ be the
semistandard tableau constructed from the successive products of
$\overline{c}_j$'s defined via the row insertion as in \cite{F} p.11.
By virtue of the $U'_q(A_{M-1})$-invariance \cite{FOY}, it is
a conserved quantity under any time evolution
$T_\kappa$, i.e.,
${\cal T}(\mbox{\mathversion{bold}$p$}) =
{\cal T}(T_\kappa(\mbox{\mathversion{bold}$p$}))$.
In the context of the Robinson-Schensted-Knuth correspondence,
${\cal T}(\mbox{\mathversion{bold}$p$})$ stands for the $P$-symbol.
For any 1 soliton state
$\mbox{\mathversion{bold}$p$} = \iota^{(L_0,L_1)}_k(b), \, b \in B'_k$,
one has ${\cal T}(\mbox{\mathversion{bold}$p$}) = b$.
One can also check that ${\cal T}(\mbox{\mathversion{bold}$p$})$ equals
\begin{center}
\setlength{\unitlength}{5mm}
\begin{picture}(2,2)(0,0)
\put(0,0){\line(1,0){1}}
\put(0,1){\line(1,0){2}}
\put(0,2){\line(1,0){2}}
\put(0,0){\line(0,1){2}}
\put(1,0){\line(0,1){2}}
\put(2,1){\line(0,1){1}}
\put(0.3,0.26){2}
\put(0.3,1.26){1}
\put(1.3,1.26){3}
\end{picture}
\hskip0.6cm
\setlength{\unitlength}{5mm}
\begin{picture}(3,2)(0,0)
\put(0,1){\line(1,0){3}}
\put(0,2){\line(1,0){3}}
\put(0,1){\line(0,1){1}}
\put(1,1){\line(0,1){1}}
\put(2,1){\line(0,1){1}}
\put(3,1){\line(0,1){1}}
\put(0.3,1.26){2}
\put(1.3,1,26){2}
\put(2.3,1.26){3}
\end{picture}
\hskip0.6cm
\setlength{\unitlength}{5mm}
\begin{picture}(5,2)(0,0)
\put(0,0){\line(1,0){1}}
\put(0,1){\line(1,0){5}}
\put(0,2){\line(1,0){5}}
\put(0,0){\line(0,1){2}}
\put(1,0){\line(0,1){2}}
\put(2,1){\line(0,1){1}}
\put(3,1){\line(0,1){1}}
\put(4,1){\line(0,1){1}}
\put(5,1){\line(0,1){1}}
\put(0.3,0.26){3}
\put(0.3,1.26){1}
\put(1.3,1.26){1}
\put(2.3,1.26){2}
\put(3.3,1.26){2}
\put(4.3,1.26){3}
\end{picture}
\end{center}
in Examples \ref{ex:tns}, \ref{ex:tnsd} and
\ref{ex:doublescattering}, respectively
throughout the scattering.

%
%
\section{$A^{(1)}_M$ Automaton as Ultradiscrete KP equation}\label{sec:toki}
Here we investigate the $A_M^{(1)}$ automaton constructed in
Section \ref{subsec:cellauto}
>from the viewpoint of ultradiscretization \cite{TTMS,MSTTT}.
With
the same notations as (\ref{eq:action}) we define $\DIS u_{n,j}^t$ and $\DIS
v_{n,j}^t$ to be the multiplicities of $(M+2-j)$th content of $\DIS b_n^t$ and
$\DIS v_n^t$, {\it i.e.},
\begin{eqnarray*}
b_n^t &=&(u_{n,M+1}^t, u_{n,M}^t, \cdots, u_{n,1}^t ),\\
v_n^t &=&(v_{n,M+1}^t, v_{n,M}^t, \cdots, v_{n,1}^t).
\end{eqnarray*}
\begin{proposition}\label{prop:BBS-rule}
The map
\begin{equation}
R \, : \, v^t_{n} \ot b^{t}_n
\stackrel{\sim}{\mapsto} b^{t+1}_n \ot v^{t}_{n+1},
\end{equation}
is expressed by   $(1 \le j \le M)$
\begin{eqnarray}
&&u_{n,j}^{t+1}-v_{n,j}^t =\max[X_1-\theta_n, X_2-\theta_n, \cdots,
X_{j-1}-\theta_n, X_j-\kappa_t, \cdots, X_M-\kappa_t,0] \nonumber \\
&&\qquad -\max[X_1-\theta_n, X_2-\theta_n, \cdots, X_j-\theta_n,
X_{j+1}-\kappa_t, \cdots, X_M-\kappa_t,0], \label{hakoEQ1}\\
&& v_{n+1,j}^t = u_{n,j}^t+v_{n,j}^t-u_{n,j}^{t+1},
\label{hakoEQ2}
\end{eqnarray}
where $\DIS X_{\ell} = X_{n;\ell}^t := \sum_{i=\ell}^M
u_{n,i}^t+\sum_{i=1}^{\ell} v_{n,i}^t$.
Actually (\ref{hakoEQ2}) is valid also for $j = M+1$.
\end{proposition}
\par\noindent
{\em Proof.}
In the present proof, we abbreviate $\DIS u_{n,j}^t$ and $\DIS v_{n,j}^t$ to
$\DIS u_j$ and $\DIS v_j$ respectively.
We also put $u_{j+M+1} = u_j, \ v_{j+M+1}=v_j$ {\it etc.},  {\it i.e.}
each suffix is defined modulo $M+1$.

We define $\DIS u_j^{(k)}, \ v_j^{(k)}$ ($j=1,2,...,M+1$, $k=1,2,...,M+1$) as
follows.
\begin{itemize}
\item[(1)]
Let $\DIS \Delta u_j = \Delta v_{j+1} := \min[u_j,v_{j+1}]$, and
$\DIS u_j^{(1)} := u_j-\Delta u_j$C$\DIS v_j^{(1)} := v_j-\Delta v_j$ for
$j=1,2,...,M+1$.
\item[(2)]
 For $\forall j$, we define $\DIS \Delta u_j^{(1)} = \Delta v_{j+2}^{(1)} :=
\min[u_j^{(1)},v_{j+2}^{(1)}]$, and
$\DIS u_j^{(2)} := u_j^{(1)}-\Delta u_j^{(1)}$, $\DIS v_j^{(2)} :=
v_j^{(1)}-\Delta v_j^{(1)}$.
\item[(3)]
Similarly to the step (2),  we recursively define
$\DIS \Delta u_j^{(\ell-1)} = \Delta v_{j+\ell}^{(\ell-1)} :=
\min[u_j^{(\ell-1)},v_{j+\ell}^{(\ell-1)}]$,
$\DIS u_j^{(\ell)} := u_j^{(\ell-1)}-\Delta u_j^{(\ell-1)}$ and $\DIS
v_j^{(\ell)} := v_j^{(\ell-1)}-\Delta v_j^{(\ell-1)}$ for $\DIS
\ell=2,3,\ldots, M+1$.
\end{itemize}
{}From Proposition \ref{prop:H-rule}, we see that $\DIS u_j^{(M+1)}$ and $\DIS
v_j^{(M+1)}$ are the numbers of {\it unconnected} dots in $(M+2-j)$ th box in the
column diagrams for $\DIS b^{t}_n$ and $\DIS v^{t}_n$ respectively.
See Example \ref{ex:NY}.
Noting that $\Delta u^{(M)}_j = \Delta v^{(M)}_j$ we have
\begin{equation}
\label{U-rel}
 u_{n,j}^{t+1}
 = v_j+u_j^{(M+1)}-v_j^{(M+1)} =
 v_j+u_j^{(M)}-v_j^{(M)}
\end{equation}
for $1 \le j \le M+1$.
The following formulae are easily shown by induction:
\begin{eqnarray}
u_j^{(\ell)}&=&\max \left[ \sum_{i=0}^{\ell-1}u_{j+i},
\sum_{i=1}^{\ell-1}u_{j+i}+v_{j+1}, \sum_{i=2}^{\ell-1}u_{j+i}+\sum_{i=1}^2
v_{j+i}, \cdots, u_{j+\ell-1}+\sum_{i=1}^{\ell-1}v_{j+i}, \sum_{i=1}^{\ell}
v_{j+i} \right] \nonumber \\
&&-v_{j+1}-\max\left[ \sum_{i=1}^{\ell-1}u_{j+i},
\sum_{i=2}^{\ell-1}u_{j+i}+v_{j+2}, \cdots,
u_{j+\ell-1}+\sum_{i=2}^{\ell-1}v_{j+i}, \sum_{i=2}^{\ell} v_{j+i} \right],
\\
v_j^{(\ell)}&=&\max \left[ \sum_{i=0}^{\ell-1}v_{j-i},
\sum_{i=1}^{\ell-1}v_{j-i}+u_{j-1}, \sum_{i=2}^{\ell-1}v_{j-i}+\sum_{i=1}^2
u_{j-i}, \cdots, v_{j-\ell+1}+\sum_{i=1}^{\ell-1}u_{j-i}, \sum_{i=1}^{\ell}
u_{j-i} \right] \nonumber \\
&&-u_{j-1}-\max\left[ \sum_{i=1}^{\ell-1}v_{j-i},
\sum_{i=2}^{\ell-1}v_{j-i}+u_{j-2}, \cdots,
v_{j-\ell+1}+\sum_{i=2}^{\ell-1}u_{j-i}, \sum_{i=2}^{\ell} u_{j-i} \right].
\end{eqnarray}
Noticing  $\DIS u_{j-1} = u_{j+M},\ v_{j+1} =  v_{j-M}$, we find
\begin{eqnarray}
&&u_j^{(M)} - v_j^{(M)} = \max \left[ \sum_{i=0}^{M}u_{j+i},
\sum_{i=1}^{M}u_{j+i}+v_{j+1}, \sum_{i=2}^{M}u_{j+i}+\sum_{i=1}^2 v_{j+i},
\cdots,
\right.\nonumber \\
&& \qquad \qquad \qquad \left. \sum_{i=M-1}^M
u_{j+i}+\sum_{i=1}^{M-1}v_{j+i}, u_{j+M}+\sum_{i=1}^{M} v_{j+i} \right]
\nonumber \\
&& \;\; -\max \left[ \sum_{i=0}^{M}v_{j-i}, \sum_{i=1}^{M}v_{j-i}+u_{j-1},
\sum_{i=2}^{M}v_{j-i}+\sum_{i=1}^2 u_{j-i}, \cdots,
\right. \nonumber \\
&& \qquad \qquad \qquad \left. \sum_{i=M-1}^M
v_{j-i}+\sum_{i=1}^{M-1}u_{j-i}, v_{j-M} + \sum_{i=1}^{M} u_{j-i} \right].
\end{eqnarray}
Subtracting $\DIS \sum_{i=0}^{j-1}u_i + \sum_{i=j+1}^{M+1} v_i$ {}from  the both
$\max[...]$ terms in the left hand side of the equation and using the
relations:
$\DIS u_{M+1} = \theta_n - \sum_{j=1}^M u_j$ and
$\DIS v_{M+1} = \kappa_t - \sum_{j=1}^M v_j$,
we get (\ref{hakoEQ1}) {}from  (\ref{U-rel}).
Noticing that the number of dots of two column diagrams are preserved in the
rule, we obtain (\ref{hakoEQ2}).
\hfill
$\square$
\par
Our goal in this section is to show that (\ref{hakoEQ1}) and
(\ref{hakoEQ2}) are ultradiscrete limits of the (one-constrained)
nonautonomous discrete KP equation (ndKP eq.):
\begin{eqnarray}
&&(b_{n+1}-c_{j+1}) \tau(t,n,j) \tau(t+1,n+1,j+1) + (c_{j+1}-a_{t+1})
\tau(t+1,n+1,j) \tau(t,n,j+1)  \nonumber \\
&& \qquad \quad + (a_{t+1}-b_{n+1}) \tau(t,n+1,j) \tau(t+1,n,j+1) = 0.
\label{ndKPeq}
\end{eqnarray}
Here $\DIS a_t,\ b_n,\ c_j$ are arbitrary complex parameters.
The ndKP eq.~(\ref{ndKPeq}), which is sometimes called the (nonautonomous)
Hirota-Miwa equation, is equivalent to the generating formulae of the KP
hierarchy \cite{S,DJKM}.
Its soliton solutions, Lax operators, Darboux transformations $etc.$ have
been investigated in \cite{WTS}.
We set $\DIS a_{t+1}=1+\delta_t$ and $\DIS b_{n+1}=1+\gamma_n$.
We also assume that $\DIS c_1=1, c_2=c_3=\cdots =c_{M+1}=0$ and
\begin{equation}
\label{reduction}
\tau(t,n,j+M+1)=\tau(t,n,j).
\end{equation}
The constraint (\ref{reduction}) is an analogue of $M$-reduction of the KP
hierarchy which restricts the space of transformation group of $\tau$
functions to the subgroup generated by $\DIS A_M^{(1)}$ \cite{DJKM}.
Let
\begin{eqnarray}
 U_{n,j}^t &:=
&\frac{\tau(t,n+1,j)\tau(t,n,j+1)}{\tau(t,n,j)\tau(t,n+1,j+1)}, \nonumber \\
V_{n,j}^t &:=
&\frac{\tau(t+1,n,j+1)\tau(t,n,j)}{\tau(t+1,n,j)\tau(t,n,j+1)}
\label{UVdef}
\end{eqnarray}
for $1 \le j \le M$.
We also introduce a small positive parameter $\varepsilon$, and put $\DIS
\delta_t=\exp\left[-\kappa_t/\varepsilon\right]$ and $\DIS
\gamma_n=\exp\left[-\theta_n/\varepsilon\right]$.
Then we have
\begin{theorem}\label{UDtheorem}
Let
\begin{eqnarray*}
u_{n,j}^t &= &\lim_{\varepsilon \to +0} \varepsilon \log U_{n,j}^t, \\
v_{n,j}^t &= &\lim_{\varepsilon \to +0} \varepsilon \log V_{n,j}^t
\end{eqnarray*}
be the ultradiscrete limits for $1 \le j \le M$, and
specify $u^t_{n,M+1}$ and $v^t_{n,M+1}$
by $\sum_{j=1}^{M+1}u^t_{n,j} = \theta_n$ and
$\sum_{j=1}^{M+1}v^t_{n,j} = \kappa_t$. Then
$\{u^t_{n,j} \}$ and $\{ v^t_{n,j} \}$ satisfy (\ref{hakoEQ1}) and (\ref{hakoEQ2}).
\end{theorem}
\par\noindent
{\em Proof}.
We use abbreviations: $\DIS \tau_j := \tau(t,n,j), \tau_j^t :=
\tau(t+1,n,j), \tau_{n,j} := \tau(t,n+1,j), \tau_{n,j}^t :=
\tau(t+1,n+1,j)$.
The ndKP eq.~(\ref{ndKPeq}) with the constraint (\ref{reduction}) is
rewritten as the following $M+1$ simultaneous equations:
\begin{eqnarray}
&&(1+\gamma_n) \tau_1 \tau_{n,2}^t -(1+\delta_t) \tau_{n,1}^t
\tau_2+(\delta_t-\gamma_n)\tau_{n,1} \tau_2^t=0, \nonumber \\
&&(1+\gamma_n) \tau_2 \tau_{n,3}^t -(1+\delta_t) \tau_{n,2}^t
\tau_3+(\delta_t-\gamma_n)\tau_{n,2} \tau_3^t=0, \nonumber \\
&& \qquad \cdots \nonumber \\
&&(1+\gamma_n) \tau_M \tau_{n,M+1}^t -(1+\delta_t) \tau_{n,M}^t
\tau_{M+1}+(\delta_t-\gamma_n)\tau_{n,M} \tau_{M+1}^t=0, \nonumber \\
&&\gamma_n \tau_{M+1} \tau_{n,1}^t - \delta_t \tau_{n,M+1}^t
\tau_1+(\delta_t-\gamma_n)\tau_{n,M+1} \tau_1^t=0.\label{simultaneouseq}
\end{eqnarray}
Defining
\begin{eqnarray*}
&x_1 := \tau_{n,1 }^t \tau_2 \tau_3 \ \cdots \tau_{M+1}, \quad
&y_1 := \tau_{n,1} \tau_2^t \tau_3 \tau_4 \cdots \tau_{M+1}, \\
&x_2 := \tau_1 \tau_{n,2 }^t \tau_3 \ \cdots \tau_{M+1}, \quad
&y_2 := \tau_1 \tau_{n,2} \tau_3^t \tau_4  \cdots \tau_{M+1}, \\
&\qquad \cdots \quad &\qquad \cdots \quad \\
&x_{M+1} := \tau_1 \tau_2 \tau_3 \ \cdots \tau_{n,M+1}^t, \quad
&y_{M+1} := \tau_1^t \tau_2 \tau_3 \tau_4 \cdots \tau_{n,M+1},\\
&\xv := (x_1,x_2,\cdots,x_{M+1})^T,
&\yv := (y_1,y_2,\cdots,y_{M+1})^T,
\end{eqnarray*}
we obtain
$$
\L \xv = (\delta_t-\gamma_n) \yv,
$$
where
$$
\L = \left(
\begin{array}{cccccc}
(1+\delta_t)&-(1+\gamma_n)&0&\cdots &0&0\\
0&(1+\delta_t)&-(1+\gamma_n)&0&\cdots&0\\
\vdots&\vdots&\ddots&\ddots&\ddots&\vdots\\
0&0&0&\cdots &(1+\delta_t)&-(1+\gamma_n)\\
-\gamma_n&0&0&\cdots&0&\delta_t
\end{array}
\right).
$$
Its inverse matrix is easily calculated as
$$
\DIS \L^{-1} = \D/\left( (1+\delta_t)^M \delta_t-(1+\gamma_n)^M \gamma_n
\right),
$$
$$
(\D)_{i,j}=\left\{
\begin{array}{lc}
(1+\gamma_n)^{M+1-i}(1+\delta_t)^{i-1}, &\quad j=M+1,\\
\delta_t(1+\delta_t)^{M+i-j-1}(1+\gamma_n)^{j-i}, &\quad j \ge i \;\;(j \ne
M+1), \\
\gamma_n(1+\gamma_n)^{M-i+j}(1+\delta_t)^{i-j-1}, &\quad j \le i-1 \;\;(j
\ne M+1).
\end{array}
\right.
$$
Thus, for $\DIS 0 < \delta_t, \gamma_n \ll 1$, we have
$$
(\delta_t-\gamma_n) \L^{-1} \sim
\left(
\begin{array}{cccccc}
\delta_t&\delta_t&\delta_t&\cdots &\delta_t &1\\
\gamma_n&\delta_t&\delta_t&\cdots&\delta_t&1\\
\gamma_n&\gamma_n&\delta_t&\cdots&\delta_t&1\\
\vdots&\vdots&\vdots&\ddots&\vdots&\vdots\\
\gamma_n&\gamma_n&\gamma_n&\cdots&\gamma_n&1
\end{array}
\right).
$$
Precisely speaking, $\DIS A \sim B$ means $\DIS \lim_{\varepsilon \to +0}
\varepsilon \log A(\varepsilon) = \lim_{\varepsilon \to +0} \varepsilon \log
B(\varepsilon)$.
Since
\begin{eqnarray}
x_j&\sim&\gamma_n \sum_{i=1}^{j-1}y_i+\delta_t\sum_{i=j}^{M}y_i +y_{M+1}, \\
x_{j+1}&\sim&\gamma_n \sum_{i=1}^{j}y_i+\delta_t\sum_{i=j+1}^{M}y_i +y_{M+1},
\end{eqnarray}
we have
\begin{equation}
\label{xequation}
\frac{x_j}{x_{j+1}} \sim \frac{\gamma_n \sum_{i=1}^{j-1} (y_i/y_{M+1})
+ \delta_t \sum_{i=j}^{M}  (y_i/y_{M+1}) \ +1}{\gamma_n \sum_{i=1}^{j}
(y_i/y_{M+1})
+ \delta_t \sum_{i=j+1}^{M}  (y_i/y_{M+1}) \ +1}.
\end{equation}
{}From the definition of $\DIS U_{n,j}^t$ and $\DIS V_{n,j}^t$, we find that
the left hand side of (\ref{xequation}) is equal to $\DIS
U_{n,j}^{t+1}/V_{n,j}^t$ and that $\DIS (y_j/y_{M+1})=\prod_{i=j}^M
U_{n,i}^t \prod_{i=1}^j V_{n,i}^t$.
Since it holds that
$$
\lim_{\varepsilon \to +0} \varepsilon \log \left(\frac{x_j}{x_{j+1}}\right)
=
\lim_{\varepsilon \to +0} \varepsilon \log \left[\mbox{\rm right hand side
of (\ref{xequation})}\right],
$$
we have (\ref{hakoEQ1}) by putting
\begin{eqnarray*}
u_{n,j}^t &= &\lim_{\varepsilon \to +0} \varepsilon \log U_{n,j}^t, \\
v_{n,j}^t &= &\lim_{\varepsilon \to +0} \varepsilon \log V_{n,j}^t.
\end{eqnarray*}
{}From the definitions (\ref{UVdef}), we have
$$
\frac{U_{n,j}^{t+1}}{U_{n,j}^{t}}=\frac{V_{n,j}^t}{V_{n+1,j}^t},
$$
which gives (\ref{hakoEQ2}) in the ultradiscrete limit.
\hfill
$\square$
\par

%
%

Next, we consider soliton solutions to the $A_M^{(1)}$ automaton.
It is obvious that if the limit :
\begin{equation}
Y_{n,j}^t :=
\lim_{\ve \to +0} \ve \log \tau(t,n,j)
\end{equation}
exists, then {}from (\ref{UVdef}) we have for $1 \le j \le M$
\begin{eqnarray}\label{eq:uvy}
u_{n,j}^t &= &Y_{n+1,j}^t+Y_{n,j+1}^t-Y_{n,j}^t-Y_{n+1,j+1}^t, \nonumber \\
v_{n,j}^t &= &Y_{n,j+1}^{t+1}+Y_{n,j}^t-Y_{n,j}^{t+1}-Y_{n,j+1}^t.
\end{eqnarray}
{}From Theorem \ref{UDtheorem}, they satisfy (\ref{hakoEQ1}) and
(\ref{hakoEQ2}).
Hence we have only to know $\DIS Y_{n,j}^t$ to get solutions to
(\ref{hakoEQ1}) and (\ref{hakoEQ2}).
We will call $\DIS Y_{n,j}^t$ an $N$ soliton solution to the $\DIS
A_M^{(1)}$ automaton when it is an ultradiscrete limit of one parameter
($\ve$) family of a certain
$M \times N$ soliton solutions $\DIS \tau(t,n,j)$ to the ndKP
eq.~(\ref{ndKPeq}) as explained in Appendix \ref{app:b}.
It indeed corresponds to an
$N$ soliton state in the sense of Section \ref{subsec:solitons}.

The following fact is well known \cite{DJKM,WTS}.
\begin{proposition}
\label{N-tau-soliton}
The $N$ soliton solution to (\ref{ndKPeq}) is given by the vacuum
expectation value:
\begin{eqnarray}
\tau(t,n,j)&=& \langle vac| g(\t)| vac \rangle, \\
g(\t)&=&\prod_{k=1}^N\left( 1+\alpha_k \psi(p_k,\t)\psi^*(q_k,\t) \right).
\label{Nsoliton}
\end{eqnarray}
Here $\t = (t,n,j)$ and
$\DIS \alpha_k$ $\DIS \ (k=1, 2, \cdots ,N)\ $ are arbitrary complex
constants.
\begin{eqnarray*}
\psi(p,\t)&=&\left[\prod_{t'}^t(a_{t'}-p)
\prod_{n'}^n(b_{n'}-p)^{-1}\prod_{j'=1}^j(-c_{j'}+p)^{-1}\right]  \psi(p),
\\
\psi^*(q,\t)&=&\left[\prod_{t'}^t(a_{t'}-q)^{-1}\prod_{n'}^n(b_{n'}-q)
\prod_{j'=1}^j(-c_{j'}+q) \right] \psi^*(q),
\end{eqnarray*}
with
$$
\prod_{n'}^n X_{n'} := \left\{
\begin{array}{cl}
\prod_{n'=1}^n X_{n'} &\quad 1 \le n\\
1 &\quad n=0\\
\prod_{n'=n+1}^0 X_{n'}^{-1} &\quad n \le -1,
\end{array}
\right.
$$
and $\DIS \psi(p),\ \psi^*(q)$ are fermionic field operators which satisfy
\begin{eqnarray*}
&&\left\{\psi(p),\psi(p')\right\}_+ := \psi(p) \psi(p')+\psi(p')
\psi(p)=0,\\
&&\left\{\psi^*(q),\psi^*(q')\right\}_+ =0, \quad
\left\{\psi(p),\psi^*(q)\right\}_+ =0 \:\:\: \mbox{\rm for $(p \ne q),$} \\
&&\langle vac| \psi(p_1) \psi(p_2) \cdots \psi(p_r)\psi^*(q_r)
\psi^*(q_{r-1}) \cdots \psi^*(q_1) |vac \rangle \\
&&= \mbox{\rm det} \left( \frac{1}{p_i-q_j} \right)_{1\le i,j \le r}
=\frac{\prod_{i<j}(p_i-p_j)(q_j-q_i)}{\prod_{i,j}(p_i-q_j)}.
\end{eqnarray*}
\end{proposition}

The $N$ soliton solution (\ref{Nsoliton}) is also a solution to
(\ref{simultaneouseq}) when it satisfies the constraint
(\ref{reduction}). We can easily show
\begin{proposition}
\label{pqconstraint}
The constraint (\ref{reduction}) is achieved if it holds that
\begin{equation}
\label{PQrelation}
\left(\frac{q_k}{p_k}\right)^M \left( \frac{1-q_k}{1-p_k} \right) = 1 \quad
(k=1,2,\cdots,N) .
\end{equation}
\end{proposition}
Note that, for a given $\DIS p_k$, there are $M$ $\DIS q_k$'s which satisfy
(\ref{PQrelation}) and $\DIS q_k \ne p_k$. We use this fact to construct
explicit solutions.

{}From Propositions \ref{N-tau-soliton} and \ref{pqconstraint}, we can
construct a class of $N$ soliton solutions to the $A_M^{(1)}$ automaton.
The result is summarized as
\begin{theorem}\label{prop:N-soliton}
\begin{equation}\label{eq:result}
Y_{n+1,j+1}^{t+1}=\max_{\vec{\mu}}\left[ \sum_{i=1}^N
\mu_i K^{(i)}(t,n,j)-A(\vec{\mu};j)\right]
\end{equation}
is an $N$ soliton solution to the $\DIS A_M^{(1)}$ automaton.
Here $\DIS \vec{\mu} = (\mu_1,\mu_2, ..., \mu_N)$ ($\DIS \mu_i=0,1$)
and  $\DIS \max_{\vec{\mu}}[ \cdots ]$ denotes the maximum among the $2^N$
values obtained by putting $\mu_i=0$ or $1$ for $i=1,2,\cdots,N$.
$$
K^{(i)}(t,n,j)
= K_0^{(i)}-\sum_{j'=1}^j \ell_{j'}^{(i)} -\sum_{t'}^t
\min[\kappa_{t'}, L^{(i)}]+ \sum_{n'}^n \min[\theta_{n'}, L^{(i)}],
$$
where the sums here are generally defined by
$$
\sum_{n'}^n X_{n'} := \left\{
\begin{array}{cl}
\sum_{n'=1}^n X_{n'} &\quad 1 \le n\\
0 &\quad n=0\\
-\sum_{n'=n+1}^0 X_{n'} &\quad n \le -1.
\end{array}
\right.
$$
$\DIS L^{(i)}, \ \ell_j^{(i)}$ $(1 \le i \le N, 1 \le j \le M)$
are non-negative integers which satisfy
$\DIS L^{(i)}=\sum_{j=1}^M \ell_j^{(i)}$,
\begin{eqnarray*}
&&L^{(1)} \ge L^{(2)} \ge \cdots \ge L^{(N)}, \\
&&\ell_j^{(1)} \ge \ell_j^{(2)} \ge \cdots \ge \ell_j^{(N)}, \quad
(j=1,2,\cdots, M),
\end{eqnarray*}
and $K_0^{(i)}$ is an arbitrary integer.
In the case: $\DIS
\left\{
\begin{array}{cl}
\mu_i=1 &\quad \mbox{\rm for} \;\; i=i_1,i_2,\cdots,i_p \\
\mu_i=0 & \quad \mbox{\rm otherwise},
\end{array}
\right.$ \\
\noindent
the phase factor $A(\vec{\mu};j)$ is given by
$$
A(\vec{\mu};j) = \sum_{k=1}^p (k-1) L^{(i_k)}
+ \sum_{k=1}^p\left( X^{(i_k)}(j+k-1)-X^{(i_k)}(j) \right),
$$
where $\DIS X^{(i)}(j) = \sum_{j'=1}^j \ell^{(i)}_{j'}$ with $\DIS
\ell_{j+M}^{(i)}=\ell_j^{(i)}$.
\end{theorem}
The proof of this theorem is parallel to that in \cite{TTM}.
We give the detail in Appendix \ref{app:b}.
For $N=1$ it is the general solution, and we conjecture that
it is also so for $M=1$.
Except these cases the above result does not cover the
arbitrary initial condition.
There is some freedom to employ different `phase factor' $A(\vec{\mu};j)$
than the  above one depending on the way in taking the ultradiscrete limit.

\section{Summary}\label{sec:summary}
In this paper we have introduced
the $A^{(1)}_M$ automaton, which is a crystal theoretic formulation of the
generalized box-ball systems.
In terms of the box-ball systems, it
corresponds to the dynamics of $M$ kinds of balls,
where the carriers and boxes have arbitrary and inhomogeneous capacities.
We have introduced the solitons labeled with the crystals $B'_k$
of $U'_q(A^{(1)}_{M-1})$.
Scattering matrices of two solitons are identified with the
combinatorial $R$ matrices of $U'_q(A^{(1)}_{M-1})$-crystals.
Piecewise linear evolution equations are obtained and
identified with an ultradiscrete limit of the nonautonomous discrete KP
equation.
It allowed us to construct a class of $N$ soliton solutions.
We have left the studies of phase shifts in the scattering
and  construction of $N$ soliton solutions corresponding to arbitrary
initial conditions for $N\ge 2$ as future problems.
The interplay between the ultradiscrete limit
of the classical integrable systems and the $q \rightarrow 0$ limit of the quantum
integrable systems elucidated in this paper deserves further investigation.

\vspace{0.4cm}
\noindent
{\bf Acknowledgements}. \hspace{0.1cm}
The authors thank M. Okado and Y. Yamada for discussions and
sending a preprint prior to the publication.
Thanks are also due to J. Matsukidaira, A. Nagai, J. Satsuma and D. Takahashi for
helpful discussions about soliton solutions.

\appendix
\section{Proof of Proposition \ref{prop:toki}}\label{app:a}

First we show that it suffices to prove Proposition \ref{prop:toki}
for  $M=1$ and $h=k$.
Without a loss of generality we may set $M=2$ and consider the
time evolution $T_{\kappa = \infty}$.
We find it convenient to adopt the equivalent box-ball system picture
explained in Section \ref{subsec:equiva}.
Thus the elements in $B^{\ot {\cal L}}_1$ in
(\ref{eq:Nsoliton}) will be represented as $\ldots 131..2 \ldots$ for example.
It stands for the array of the balls
with the indices 1,3,1 and 2
and $.$ denotes an empty box.
(So they do {\em not} correspond to the letters in the
semistandard tableaux in the crystal notation.)
We keep the same notation ${\hat \theta}$ to denote the map
corresponding to (\ref{eq:thetahat}) in the box-ball picture.
It groups the array of balls and empty boxes locally together into the boxes
with capacities $\ldots, \theta_n, \theta_{n+1},\ldots$.
Then the assertion of Proposition \ref{prop:toki}
is that the scattering
\begin{equation}\label{eq:sca1}
{\hat \theta}(\ldots\overbrace{2 \cd 2}^l .\ldots
\overbrace{1 \cd 1}^{k-h}\overbrace{2 \cd 2}^h\ldots )
\xrightarrow{({\tilde T}_2{\tilde T}_1)^t}
{\hat \theta}(\ldots\overbrace{2 \cd 2}^k .\ldots
\overbrace{1 \cd 1}^{k-h}\overbrace{2 \cd 2}^{l+h-k}\ldots )
\end{equation}
takes place for sufficiently large $t$.
Here ${\tilde T}_1, {\tilde T}_2$ are the ball-moving
operators defined in Section \ref{subsec:equiva}, and
we have used $T^t_\infty = ({\tilde T}_2{\tilde T}_1)^t$  in view of
(\ref{eq:factor}) and the fact that the balls with index $\ge 3$ are absent.
In (\ref{eq:sca1}) the sequences $.\ldots$ of the empty boxes are sufficiently
long since both sides are to represent the asymptotic 2 soliton states
in the sense of Section \ref{subsec:solitons}.
Now we make use of the relation
$({\tilde T}_2{\tilde T}_1)^t =
{\tilde T}_2({\tilde T}_1{\tilde T}_2)^{t-1}
{\tilde T}_1$.
{}From the definition of the operators ${\tilde T}_i$'s and
the assumption that the 2 solitons are enough separated, (\ref{eq:sca1})
is equivalent to
\begin{equation}\label{eq:sca2}
{\hat \theta}(\ldots\overbrace{2 \cd 2}^l .\ldots
\overbrace{2 \cd 2}^h\overbrace{1 \cd 1}^{k-h}\ldots )
\xrightarrow{({\tilde T}_1{\tilde T}_2)^{t-1}}
{\hat \theta}(\ldots\overbrace{2 \cd 2}^k .\ldots
\overbrace{2 \cd 2}^{l+h-k}\overbrace{1 \cd 1}^{k-h}\ldots ).
\end{equation}
But this is justified once one establishes
\begin{equation}\label{eq:sca3}
{\hat \theta}(\ldots\overbrace{1 \cd 1}^l .\ldots
\overbrace{1 \cd 1}^k\ldots )
\xrightarrow{{\tilde T}^{t-1}_1}
{\hat \theta}(\ldots\overbrace{1 \cd 1}^k .\ldots
\overbrace{1 \cd 1}^{l}\ldots ),
\end{equation}
because (\ref{eq:sca2}) and (\ref{eq:sca3})
correspond to the same canonical system
\par\noindent
${\hat \theta}(\ldots {12 \cd l} .\ldots
{l+1 \cd l+k}\ldots )$
in the sense of
Section \ref{subsec:equiva} with respect to the relevant
time evolutions and therefore they possess
the parallel time evolution pattern owing to (\ref{eq:canonicalcom}).
In this way
the proof of Proposition \ref{prop:toki} is reduced to
(\ref{eq:sca3}), which is equivalent to the case $M=1$ and $h=k$.

Now setting  $L^{(1)}= l$ and $L^{(2)}=k$,  we are to show
\begin{proposition}\label{prop:M1proof}
Set  $M=1$, assume that $\DIS \kappa_t \gg L^{(1)}$ ($\forall t$)
and $\DIS \theta_n < L^{(1)}$ for all but finitely many $n$'s.
Then two solitons with amplitudes $\DIS L^{(1)}$ and
$\DIS L^{(2)}$ ($\DIS L^{(1)} > L^{(2)}$)
scatter into two solitons with amplitudes
$\DIS L^{(2)}$ and $\DIS L^{(1)}$, respectively.
\end{proposition}
Namely, the amplitudes of two solitons do not change after the collision.
To prove the proposition, we need several lemmas.
The following two lemmas are obvious.
\begin{lemma}\label{M1lemma1}
For given integers $K_1$ and $K_2$, if there exists an integer $n_0$ such that
$$
K_1 + \sum_{n'}^{n_0} \min\left[\theta_{n'}, L^{(1)} \right]
\ge 0 \ge
K_2 + \sum_{n'}^{n_0} \min\left[\theta_{n'}, L^{(2)} \right],
$$
then, for $n \ge n_0$,
$$
K_1 + \sum_{n'}^{n} \min\left[\theta_{n'}, L^{(1)} \right]
\ge
K_2 + \sum_{n'}^{n} \min\left[\theta_{n'}, L^{(2)} \right],
$$
and for $ n < n_0$
$$
0 > K_2 + \sum_{n'}^{n} \min\left[\theta_{n'}, L^{(2)} \right].
$$
\end{lemma}
\begin{lemma}\label{M1lemma2}
For given integers ${K'}_1$ and ${K'}_2$, if there exists an integer $n_0$ such that
$$
{K'}_2 + \sum_{n'}^{n_0} \min\left[\theta_{n'}, L^{(2)} \right] \ge 2L^{(2)}> 0 \ge {K'}_1 + \sum_{n'}^{n_0} \min\left[\theta_{n'}, L^{(1)} \right],
$$
then, for $n \ge n_0$,
$$
{K'}_1+{K'}_2 +\sum_{n'}^{n} \min\left[\theta_{n'}, L^{(1)} \right]
+ \sum_{n'}^{n} \min\left[\theta_{n'}, L^{(2)} \right]-2L^{(2)}
\ge
{K'}_1 + \sum_{n'}^{n} \min\left[\theta_{n'}, L^{(1)} \right],
$$
and for $n < n_0$,
$$
0 > {K'}_1 + \sum_{n'}^{n} \min\left[\theta_{n'}, L^{(1)} \right].
$$
\end{lemma}

Now we define an integer $\DIS N_0(t)$ for given integers $\DIS K_2$ and $t$ as
\begin{eqnarray*}
&&K_2-L^{(2)} t + \sum_{n'}^{N_0(t)}\min\left[\theta_{n'}, L^{(2)} \right] \\
&&\ge 2 L^{(2)} \\
&&> K_2-L^{(2)} t + \sum_{n'}^{N_0(t)-1}\min\left[\theta_{n'}, L^{(2)} \right].
\end{eqnarray*}
With this $\DIS N_0(t)$ we can show
\begin{lemma}\label{M1lemma3}
For any  integers $\DIS K_1$ and $\DIS K_2$, we have
\begin{align*}
\lim_{T \rightarrow \infty}& (
K_2 -L^{(2)}T + \sum_{n'}^{N_0(T)}\min\left[\theta_{n'}, L^{(2)} \right]  \\
&-K_1 +L^{(1)}T - \sum_{n'}^{N_0(T)}\min\left[\theta_{n'}, L^{(1)} \right] )
= +\infty.
\end{align*}
\end{lemma}
\par\noindent
{\em Proof}.
{}From the definition of $\DIS N_0(t)$,  we have
\begin{equation}
\label{M1L3-1}
-L^{(2)} < -t L^{(2)} + \sum_{n'=N_0(0)+1}^{N_0(t)}\min\left[\theta_{n'},L^{(2)}\right] < L^{(2)}.
\end{equation}
Hence we have
\begin{eqnarray}
\Delta(t) &:= &-L^{(2)}t + \sum_{n'=N_0(0)+1}^{N_0(t)}\min\left[\theta_{n'}, L^{(2)} \right] \nonumber\\
&&\qquad -\left( -L^{(1)}t + \sum_{n'=N_0(0)+1}^{N_0(t)}\min\left[\theta_{n'}, L^{(1)} \right] \right) \nonumber\\
&= &t\left(L^{(1)}-L^{(2)}\right)-\sum_{n'=N_0(0)+1}^{N_0(t)}\left(\min\left[\theta_{n'}, L^{(1)} \right] - \min\left[\theta_{n'}, L^{(2)} \right] \right) \nonumber\\
&\ge &t\left(L^{(1)}-L^{(2)}\right)-\sum_{n'=N_0(0)+1}^{N_0(t)}\left(\theta_{n'} - \min\left[\theta_{n'}, L^{(2)} \right] \right) \nonumber\\
&> &t L^{(1)}-L^{(2)}-\sum_{n'=N_0(0)+1}^{N_0(t)}\theta_{n'}.
\label{M1L3-2}
\end{eqnarray}
{}From (\ref{M1L3-1}), we obtain an inequality:
$$
t> -1+ \frac{1}{L^{(2)}} \sum_{n'=N_0(0)+1}^{N_0(t)}\min\left[\theta_{n'}, L^{(2)} \right].
$$
Thus, {}from (\ref{M1L3-2}), we find
\begin{eqnarray}\label{eq:almost}
\Delta(t) &>&
-(L^{(1)}+L^{(2)}) + \sum_{n'=N_0(0)+1}^{N_0(t)}\min\left[\frac{L^{(1)}-L^{(2)}}{L^{(2)}}\theta_{n'}, L^{(1)}-\theta_{n'} \right].
\end{eqnarray}
Since $\DIS L^{(1)}>L^{(2)}$,  $\DIS L^{(1)} > \theta_n$ for
all but finitely many $n$'s
and $\DIS \lim_{t \to +\infty} N_0(t) = +\infty$
which is seen {}from (\ref{M1L3-1}), we find
\begin{equation}
\label{M1L3-3}
\lim_{t \to +\infty} \Delta(t) = +\infty.
\end{equation}
This suffices to prove the lemma.
\hfill
$\square$
\par
Now we prove Proposition \ref{prop:M1proof}.
{}From  (\ref{eq:uvy}) we have
\begin{eqnarray}
u_n^t&:= &u_{n,j=1}^{t+1}  \nonumber \\
& = & Y_{n+1,1}^{t+1}-{Y}_{n+1,2}^{t+1}-Y_{n,1}^{t+1}+{Y}_{n,2}^{t+1}.\label{M12sol}
\end{eqnarray}
Specializing Theorem \ref{prop:N-soliton} to
a two soliton solution with
$M=1$ and $\kappa_t=+\infty$,  we have
\begin{eqnarray*}
Y_{n+1,1}^{t+1}
&= & \max \left[0, K_1(n,t), K_2(n,t), K_1(n,t)+K_2(n,t)-2L^{(2)} \right], \\
{Y}_{n+1,2}^{t+1}
&= & \max \left[0, K_1(n,t)-L^{(1)}, K_2(n,t)-L^{(2)}, K_1(n,t)+K_2(n,t)-L^{(1)}-3L^{(2)} \right], \\
K_{i}(n,t) &= & K_{i} +\sum_{n'}^{n} \min\left[ \theta_{n'}, L^{(i)} \right]-t L^{(i)} \quad (i=1,2).
\end{eqnarray*}
Note that $Y^t_{n,2} = Y^{t+1}_{n,1}$  due to the last equation in
(\ref{simultaneouseq}) and the condition $\delta_t =\exp[-\kappa_t/\varepsilon]=0$.
Given $L^{(1)} > L^{(2)}$,
there exist integers $n_1, \ n_2,\ j, \ r_1, \ r_2$ that  satisfy
\begin{align*}
&n_1 \ll n_1 + j \ll n_2, \quad 1 \le r_1 \le \min(L^{(1)}, \theta_{n_1}),
\quad 1 \le r_2 \le \min(L^{(2)},\theta_{n_2}),\\
&K_1 + \sum_{n'}^{n_1} \min\left[ \theta_{n'}, L^{(1)}\right] > 0
\ge K_2 + \sum_{n'}^{n_1} \min\left[ \theta_{n'}, L^{(2)}\right],\\
&K_1-L^{(1)} + \sum_{n'}^{n_1+j} \min\left[ \theta_{n'}, L^{(1)}\right] > 0
\ge K_2-L^{(2)} + \sum_{n'}^{n_1+j} \min\left[ \theta_{n'}, L^{(2)} \right],
\end{align*}
where $K_i$ $(i=1,2)$ is defined by
\begin{eqnarray*}
K_1 &= & r_1-\sum_{n'}^{n_1} \min\left[ \theta_{n'}, L^{(1)}\right], \\
K_2 &= & r_2+2L^{(2)} -\sum_{n'}^{n_2} \min\left[ \theta_{n'}, L^{(2)}\right].
\end{eqnarray*}
{}From Lemma \ref{M1lemma1}, we find at $t=0$ that
\begin{eqnarray}
Y_{n+1,1}^{1} &= & \max \left[0, K_1(n,0), K_1(n,0)+K_2(n,0)-2L^{(2)} \right], \nonumber \\
{Y}_{n+1,2}^{1} &= & \max \left[0, K_1(n,0)-L^{(1)}, K_1(n,0)+K_2(n,0)-L^{(1)}-3L^{(2)} \right]. \label{M1Yval}
\end{eqnarray}
Substituting (\ref{M1Yval}) into eq.~(\ref{M12sol}), we obtain
\begin{equation}
u_n^0 =
\left\{
\begin{array}{lcl}
0 & \quad & \mbox{\rm for $n < n_1$}\\
r_1 & \quad & \mbox{\rm for $n = n_1$}\\
\theta_n & \quad & \mbox{\rm for $n_1< n < {n'}_1$}\\
L^{(1)}-\sum_{n'=n_1+1}^{{n'}_1-1}\theta_{n'}-r_1 & \quad & \mbox{\rm for $n = {n'}_1$}\\
0 & \quad & \mbox{\rm for ${n'}_1<n < n_2$}\\
r_2 & \quad & \mbox{\rm for $n = n_2$}\\
\theta_n & \quad & \mbox{\rm for $n_2< n < {n'}_2$}\\
L^{(2)}-\sum_{n'=n_2+1}^{{n'}_2-1}\theta_{n'}-r_2 & \quad & \mbox{\rm for $n = {n'}_2$}\\
0 & \quad & \mbox{\rm for ${n'}_2<n $},
\end{array}
\right.
\end{equation}
where ${n'}_i$ ($i=1,2$) are defined by $n'_i = n_i+1$ if
$r_i = L^{(i)}$, and otherwise by
$$
L^{(i)}-\sum_{n'=n_i+1}^{{n'}_i}\theta_{n'}-r_i  \le 0 <L^{(i)}-\sum_{n'=n_i+1}^{{n'}_i-1}\theta_{n'}-r_i.
$$
Thus we see that the two soliton solution can correspond to $any$ initial configuration in which $L^{(1)}$ soliton is situated left hand side of $L^{(2)}$ soliton with sufficient spacing.
Hence, to prove the proposition, we have only to show that the solution $\DIS u_n^t$ describes the two soliton state in which $L^{(2)}$ soliton is left hand side of $L^{(1)}$ soliton for $t \gg 1$.
\par
{}From the definition of $\DIS N_0(t)$ and Lemma \ref{M1lemma3},
there exists $T$ and $j$ such that
\begin{align*}
&K_2 -L^{(2)}T + \sum_{n'}^{N_0(T)}\min\left[\theta_{n'}, L^{(2)} \right]
> 0  \ge
K_1 -L^{(1)}T + \sum_{n'}^{N_0(T)}\min\left[\theta_{n'}, L^{(1)} \right],\\
&K_2 -L^{(2)}T-L^{(2)} + \sum_{n'}^{N_0(T)+j}\min\left[\theta_{n'}, L^{(2)} \right]
> 0  \ge
K_1 -L^{(1)}T-L^{(1)} + \sum_{n'}^{N_0(T)+j}\min\left[\theta_{n'}, L^{(1)} \right].
\end{align*}
Thus, {}from Lemma \ref{M1lemma2}, we have at $t = T$ that
\begin{eqnarray*}
Y_{n+1,1}^{T+1} &= & \max \left[0, K_2(n,T), K_1(n,T)+K_2(n,T)-2L^{(2)} \right], \\
{Y}_{n+1,2}^{T+1} &= & \max \left[0, K_2(n,T)-L^{(2)}, K_1(n,T)+K_2(n,T)-L^{(1)}-3L^{(2)} \right].
\end{eqnarray*}
Substituting these into eq.~(\ref{M12sol}), we find that $\DIS u_n^T$
describes a configuration in which $L^{(2)}$ soliton locates around
$n=N_0(T)$ and $L^{(1)}$ soliton does around $n \gg N_0(T)$.
This completes the proof.

\section{Derivation of $N$ soliton solutions}\label{app:b}

Here we explain the derivation of the
$N$ soliton solution in Theorem  \ref{prop:N-soliton}
along the simple cases $N=1$ and $N=2$.
First we consider one soliton solution.
We will show that it has the form:
\begin{equation}
\label{1-sol.sol}
Y_{n,j}^{t}=\max\left[0,\  K_0-\sum_{i=1}^{j-1} \ell_i-\sum_{t'}^{t-1}
\min[\kappa_{t'},L]+\sum_{n'}^{n-1} \min[\theta_{n'},L] \right],
\end{equation}
where $L$ is the amplitude, $\DIS K_0$ is an integer which is
related to the phase of the soliton, and $\ell_i$ ($i=1,2,\cdots,M$) are the
non-negative integers which correspond to the number of $i$th balls in the
soliton and $\DIS \sum_{i=1}^M \ell_i=L$.
We give some details of its derivation, because similar technical
difficulties in obtaining multi-soliton solutions are resolved in the same
way.

To obtain (\ref{1-sol.sol}), we take $g(\t)$ in (\ref{Nsoliton}) as
\begin{eqnarray}
g(\t)&=&\prod_{\ell=0}^{M-1} \left(
1+c_{\ell}(p)\psi(p,\t)\psi^*(q_{\ell},\t) \right)\nonumber \\
 &= &1+ \psi(p,\t)\phi^*(p,\t), \label{1soliton} \\
\phi^*(p,\t) &:= & \sum_{\ell=0}^{M-1}
c_{\ell}(p)\psi^*(q_{\ell},\t), \nonumber
\end{eqnarray}
where $\DIS q_{\ell}$ $\ (\ell=0,1,\cdots,M-1)$  are the roots of the algebraic
equation :
\begin{equation}
\frac{x^M(1-x)-p^M(1-p)}{x-p}=0, \;\; (x \ne p)
\label{qequation}
\end{equation}
for a given real number $p$ [$\DIS (1+M^{-1})^{-1} < p < 1$], and $\DIS
c_{\ell}(p)$ $(0 \le \ell \le M-1)$ are complex coefficients which will be
determined later. Since (\ref{qequation}) has one real positive root, we
assume that $\DIS q_0$ is positive and we put $\DIS \eta=q_0/p$. Then $p$
and $\DIS q_0$ satisfy
\begin{eqnarray}
p&=&\frac{1-\eta^M}{1-\eta^{M+1}}, \\
1-p&=&\eta^M\left(\frac{1-\eta}{1-\eta^{M+1}}\right), \\
q_0&=&\eta \left( \frac{1-\eta^M}{1-\eta^{M+1}} \right).
\end{eqnarray}
The $\tau$-function $\DIS \tau(t,n,j)$ is given by vacuum expectation value as
\begin{eqnarray}
\tau(t,n,j) &=&\langle vac| g(\t)| vac \rangle  \nonumber \\
&=& 1 +
\sum_{\ell=0}^{M-1}c_{\ell}(p)\frac{1}{p-q_{\ell}}\left(\frac{q_{\ell}}{p}
\right)^{j-M-1}  \left(
\frac{1-p/(1+\delta_{0})}{1-q_{\ell}/(1+\delta_{0})} \right)
\left(
\frac{1-q_{\ell}/(1+\gamma_{0})}{1-p/(1+\gamma_{0})} \right)
\nonumber \\
&& \qquad \times \prod_{t'}^{t-1}  \left(
\frac{1-p/(1+\delta_{t'})}{1-q_{\ell}/(1+\delta_{t'})} \right)
\prod_{n'}^{n-1}  \left(
\frac{1-q_{\ell}/(1+\gamma_{n'})}{1-p/(1+\gamma_{n'})} \right).
\label{sigmaformula}
\end{eqnarray}
We introduce a small positive parameter $\varepsilon$ and put
$\eta=\exp[-L/(M \varepsilon)]$. We also put
\begin{eqnarray}
\tilde{c}_{\ell}(p) &:= &\frac{c_{\ell}(p)}{p-q_{\ell}}
\left(\frac{q_{\ell}}{p}
\right)^{-M-1}\left(
\frac{1-p/(1+\delta_{0})}{1-q_{\ell}/(1+\delta_{0})} \right)
\left(
\frac{1-q_{\ell}/(1+\gamma_{0})}{1-p/(1+\gamma_{0})} \right)
\nonumber \\
&&\quad \times
\prod_{t'=-T_0}^{0} \left(1-q_{\ell}/(1+\delta_{t'}) \right)
\prod_{n'=1}^{N_0}
\left(1-q_{\ell}/(1+\gamma_{n'})\right), \\
\chi_p(s)  &:= & \sum_{\ell=0}^{M-1}
\tilde{c}_{\ell}(p)\left(\frac{q_{\ell}}{p} \right)^s \quad  ( s  \in  \Z ),
\end{eqnarray}
where $\DIS T_0=T_0(\varepsilon)$ and $\DIS N_0=N_0(\varepsilon)$ are
positive integers which satisfy $\DIS T_0 \simeq N_0 \simeq 1/\varepsilon$.
Hence, $\DIS \lim_{\varepsilon \to +0} T_0=\lim_{\varepsilon \to +0}
N_0=+\infty$.
Since
\begin{eqnarray*}
\chi_p(s+M)&=&\sum_{\ell=0}^{M-1}
\tilde{c}_{\ell}(p)\left(\frac{q_{\ell}}{p} \right)^{s+M} \\
&=&\sum_{\ell=0}^{M-1} \tilde{c}_{\ell}(p)\left(\frac{q_{\ell}}{p}
\right)^{s}\left( \frac{1-p}{1-q_{\ell}}\right) \\
&=&(1-p)\sum_{i=0}^{\infty} p^i \sum_{\ell=0}^{M-1}
\tilde{c}_{\ell}(p)\left(\frac{q_{\ell}}{p} \right)^{s+i} \\
&=&(1-p)\sum_{i=0}^{\infty} p^i\chi_p(s+i),
\end{eqnarray*}
we have
\begin{equation}
\label{chiEQ}
\chi_p(s+M)=\sum_{i=0}^{M-1}\left(
\sum_{\ell=0}^{\infty}(1-p)^{\ell+1}p^{M\ell} \varrho_{\ell}(i) \right)
p^i\chi_p(s+i),
\end{equation}
where $\DIS \varrho_0 (i)=1$, $\DIS \varrho_1 (i)=i+1$ and
\begin{eqnarray*}
\varrho_{\ell}(i)&=&\sum_{k_1=(\ell-1)M}^{(\ell-1)M+i}
\sum_{k_2=(\ell-2)M}^{k_1} \cdots \sum_{k_{\ell}=0}^{k_{\ell-1}} 1\\
&=&\frac{(i+1)}{\ell !}\prod_{j=1}^{\ell-1} \left( \ell M + i+ j+1 \right),
\end{eqnarray*}
for $\DIS \ell \ge 2$.
Note that $\DIS \chi_p(s)$ is a real function when $\DIS \chi_p(j)$ $(0 \le
j \le M-1)$ are real.
The ratio $\DIS \varrho_{\ell+1}(i)/\varrho_{\ell}(i)$ ($\ell \ge 1, 0 \le i \le M-1$) is
calculated as
\begin{eqnarray*}
\frac{\varrho_{\ell+1}(i)}{\varrho_{\ell}(i)}&=&\frac{(\ell+1)(M+1)+i}{\ell+1}
\prod_{k=2}^{\ell} \left(1+\frac{M}{\ell M+i+k} \right) \\
&< &(M+1)\left(1+\frac{1}{\ell} \right)^{\ell} \\
&< &(M+1) e .
\end{eqnarray*}
Hence, if it holds that
$\DIS (1-p)p^M < (M+1)^{-1} e^{-1} $, we obtain
\begin{equation}
\label{chiINEQ}
|\chi_p(s+M)| \le (1-p)\sum_{i=0}^{M-1}\left(
1+(i+1)\frac{(1-p)p^M}{1-(1-p)p^M(M+1)e} \right)|\chi_p(s+i)|.
\end{equation}
Thus we find $\DIS \chi_p(s+M) \sim \eta^M \sum_{i=0}^{M-1}\chi_p(s+i)$ for
sufficiently small $\eta$.

We assume the following for $\DIS \chi_p(j)$:
\begin{eqnarray}
\chi_p(1)&=&\chi_0 \nonumber \\
\chi_p(2)&=&N_1 y^{\ell_1} \chi_p(1) \nonumber \\
\chi_p(3)&=&N_2 y^{\ell_2} \chi_p(2) \nonumber \\
&\cdots& \nonumber \\
\chi_p(M)&=&N_{M-1}y^{\ell_{M-1}} \chi_p(M-1).
\label{chiVAL}
\end{eqnarray}
Here $\DIS \chi_0$ is a positive number which is related to the initial
phase of soliton, $y=\exp[-1/\varepsilon]$, $\DIS \ell_j$ and $\DIS
N_j=N_j(\ve)$ $(j=1,2,\cdots,M-1)$ are non-negative integers and positive
numbers respectively.
They are also supposed to satisfy
\begin{eqnarray}
&&\ell_M := L-\sum_{j=1}^{M-1} \ell_j \ge 0, \nonumber \\
&& \lim_{\ve \to 0} \ve \log N_j(\ve) = 0, \label{Ncondition} \\
&&N_j y^{\ell_j} \le \varepsilon^{N^*}, \nonumber
\end{eqnarray}
for a sufficiently large positive integer $\DIS N^*$. {}From these conditions,
$\DIS \tilde{c}_{\ell}(p)$ $(0 \le \ell \le M-1)$ are uniquely determined by
the equation:
\begin{equation}
\left(
\begin{array}{cccc}
q_0&q_1&\cdots&q_{M-1}\\
q_0^2&q_1^2&\cdots&q_{M-1}^2\\
\vdots&\vdots&\ddots&\vdots\\
q_0^{M}&q_1^{M}&\cdots&q_{M-1}^{M}
\end{array}
\right)
\left(
\begin{array}{c}
\tilde{c}_0(p)\\
\tilde{c}_1(p)\\
\vdots\\
\tilde{c}_{M-1}(p)
\end{array}
\right)
=
\left(
\begin{array}{c}
p\chi_p(1)\\
p^2 \chi_p(2)\\
\vdots \\
p^{M} \chi_p(M)
\end{array}
\right).
\end{equation}
Note that the determinant of the $M \times M$ matrix in the left hand side
is equal to $\DIS \left( \prod_{i=0}^{M-1} q_i\right)
\left(\prod_{j>i}(q_j-q_i)\right) \ne 0$. It should be also noted {}from
(\ref{chiINEQ})--(\ref{Ncondition}) that
\begin{eqnarray}
 \chi_p(i) &\ge& \ve ^{-N^*} \chi_p(i+1) \quad \mbox{\rm for $\DIS \forall
i$}, \nonumber \\
\chi_p(i)&\ge& C \exp[L/\ve] \chi_p(i+M) \quad \mbox{\rm for $\DIS \forall
i$ and $\exists C>0$}.
\label{chiIEQ}
\end{eqnarray}
{}From (\ref{sigmaformula}), we have
\begin{eqnarray}
\tau(t,n,j) &= &1+\sum_{\ell=0}^{M-1} \tilde{c}_{\ell}(p)
\left(\frac{q_{\ell}}{p}\right)^{j} \prod_{t'=-T_0}^{t-1}
\left(1-\frac{q_{\ell}}{1+\delta_{t'}}\right)^{-1}
\prod_{t'}^{t-1}\left(1-\frac{p}{1+\delta_{t'}}\right) \nonumber \\
&& \qquad \times \prod_{n'}^{n-1}\left(1-\frac{p}{1+\gamma_{n'}}\right)^{-1}
\prod_{n'=n}^{N_0}\left(1-\frac{q_{\ell}}{1+\gamma_{n'}}\right)^{-1}.
\label{sigmaformula2}
\end{eqnarray}
Hereafter we restrict ourselves to the region: $\DIS |n| \le N_0$ and $\DIS
|t| \le T_0$.
Noticing  that
\begin{eqnarray*}
&& \prod_{t'=-T_0}^{t-1} \left(1-\frac{q_{\ell}}{1+\delta_{t'}}\right)^{-1}
\prod_{n'=n}^{N_0}\left(1-\frac{q_{\ell}}{1+\gamma_{n'}}\right)^{-1}\\
&&\quad =  1+\left(\sum_{t=-T_0}^{t-1}\left(\frac{1}{1+\delta_{t'}}\right)
+\sum_{n'=n}^{N_0}\left(\frac{1}{1+\gamma_{n'}}\right) \right)
q_{\ell}+\cdots \\
&&\quad =: 1 + a_1 \left(\frac{q_{\ell}}{p}\right)+a_2
\left(\frac{q_{\ell}}{p}\right)^2
+ a_3 \left(\frac{q_{\ell}}{p}\right)^3+\cdots,
\end{eqnarray*}
we find
\begin{equation}
\tau(t,n,j) = 1+\prod_{t'}^{t-1}\left(1-\frac{p}{1+\delta_{t'}}\right)
\prod_{n'}^{n-1} \left(1-\frac{p}{1+\gamma_{n'}}\right)^{-1}
\sum_{i=0}^{\infty} a_i \chi_p(j+i),
\end{equation}
where $\DIS a_0=1$ and $\DIS a_{i+1}/a_i \sim \varepsilon^{-1}$.
{}From (\ref{chiIEQ}), we have $\DIS 0 < \sum_{i=1}^{\infty} a_i \chi_p(j+i) <
\chi_p(j)$ for sufficiently small $\ve$. Putting $\DIS \chi_0 =
\exp\left[K_0/\ve\right]$ and noticing the relation:
\begin{eqnarray*}
&&\lim_{\ve \to +0} \ve \log (1-p)= -L, \\
&&\lim_{\ve \to +0} \ve \log \chi_p(j)=K_0-\sum_{i=1}^{j-1} \ell_i, \\
&&\lim_{\ve \to +0} \ve \log \left(1-\frac{p}{1+\gamma_{n}}\right)^{-1}
= \min\left[ L, \theta_n \right],  \\
&&\lim_{\ve \to +0} \ve \log \left(1-\frac{p}{1+\delta_{t}}\right)
= -\min\left[ L, \kappa_t \right],
\end{eqnarray*}
we obtain
\begin{equation}
\lim_{\ve \to +0} \ve \log \tau(t,n,j)=
\max\left[0, K_0-\sum_{i=1}^{j-1} \ell_i -\sum_{t'}^{t-1}
\min[\kappa_{t'},L]  +\sum_{n'}^{n-1} \min[\theta_{n'},L] \right].
\end{equation}
Since $\DIS \lim_{\ve \to +0} N_0(\ve)=\lim_{\ve \to +0} T_0(\ve)=+\infty$,
we have shown that (\ref{1-sol.sol}) is a one soliton solution to the
$A^{(1)}_M$ automaton.

%
%
Next  we consider two soliton solutions.
{}From the above arguments about one soliton solution, we see that the field
operators $\psi(p)$ and $\phi^*(p)$ are essentially determined by $\DIS L$,
$\DIS \ell_j \ (j=1,2,\cdots,M)$ and $\DIS K_0$. Therefore we denote these
operators by
\begin{equation}
\psi(p)=\psi(L:\ve), \quad \phi^*(p)=\phi^*(L; \{\ell_j\}; K_0: \ve).
\end{equation}
Then we take
\begin{equation}
g(\t)=(1+\psi(p_1,\t)\phi^*(p_1,\t))(1+\psi(p_2,\t)\phi^*(p_2,\t)),
\end{equation}
where
\begin{equation}
\psi(p_i)=\psi(L^{(i)}:\ve), \quad \phi^*(p_i)=\phi^*(L^{(i)};
\{\ell_j^{(i)}\}; K_0^{(i)}: \ve) \quad (i=1,2).
\end{equation}
We also assume $\DIS L^{(1)} \ge L^{(2)}$ and $\ell_j^{(1)} \ge
\ell_j^{(2)}$ $\ (j=1,2,\cdots,M)$. As we shall see below, the latter
condition turns out to be a natural constraint for soliton solutions.
Using the similar notations as above, we have
\begin{eqnarray}
\tau(t,n,j) &=& \langle
vac|(1+\psi(p_1,\t)\phi^*(p_1,\t))(1+\psi(p_2,\t)\phi^*(p_2,\t))|vac\rangle
\nonumber \\
&=&1+\langle vac|\psi(p_1,\t)\phi^*(p_1,\t)|vac\rangle +\langle
vac|\psi(p_2,\t)\phi^*(p_2,\t)|vac\rangle \nonumber \\
&&\quad +\langle
vac|\psi(p_1,\t)\phi^*(p_1,\t)\psi(p_2,\t)\phi^*(p_2,\t)|vac\rangle .
\end{eqnarray}
The second and third terms are calculated in the same way as above.
The fourth term is evaluated as
\begin{eqnarray}
&&\langle
vac|\psi(p_1,\t)\phi^*(p_1,\t)\psi(p_2,\t)\phi^*(p_2,\t)|vac\rangle
\nonumber \\
&&\;   = \sum_{\ell_1=0}^{M-1}\sum_{\ell_2=0}^{M-1} \tilde{c}_{\ell_1}(p_1)
\tilde{c}_{\ell_2}(p_2) \left(  \frac{(p_1-p_2)(q_{\ell_2}^{(2)}
-q_{\ell_1}^{(1)})}
{(p_1- q_{\ell_2}^{(2)})(p_2-q_{\ell_1}^{(1)})}  \right) \nonumber \\
&& \quad \times \prod_{i=1,2}
\left(\frac{q_{\ell_i}^{(i)}}{p_i}\right)^{j} \prod_{t'=-T_0}^{t-1}
\left(1-\frac{q_{\ell_i}^{(i)}}{1+\delta_{t'}}\right)^{-1}
\prod_{t'}^{t-1}\left(1-\frac{p_i}{1+\delta_{t'}}\right) \nonumber \\
&& \qquad \times
\prod_{n'}^{n-1}\left(1-\frac{p_i}{1+\gamma_{n'}}\right)^{-1}
\prod_{n'=n}^{N_0}\left(1-\frac{q_{\ell_i}^{(i)}}{1+\gamma_{n'}}\right)^{-1}
..\label{4term}
\end{eqnarray}
We define $\DIS \chi_{p_i}(s)$ by
\begin{equation}
\chi_{p_i}(s) := \sum_{\ell=0}^{M-1}
\tilde{c}_{\ell}(p_i)\left(\frac{q_{\ell}^{(i)}}{p_i} \right)^s \quad
(i=1,2),
\end{equation}
and suppose
\begin{eqnarray}
\chi_{p_i}(1)&=&\chi_0^{(i)} \nonumber \\
\chi_{p_i}(2)&=&N_1^{(i)} y^{\ell_1^{(i)}} \chi_{p_i}(1) \nonumber \\
\chi_{p_i}(3)&=&N_2^{(i)} y^{\ell_2^{(i)}} \chi_{p_i}(2) \nonumber \\
&\cdots& \nonumber \\
\chi_{p_i}(M)&=&N_{M-1}^{(i)}y^{\ell_{M-1}^{(i)}} \chi_{p_i}(M-1),
\label{chiVAL2}
\end{eqnarray}
where positive numbers $\DIS N_j^{(i)}$ satisfy the similar inequalities to
(\ref{Ncondition}).
{}From the assumption: $\DIS \ell_j^{(1)} \ge \ell_j^{(2)}$
$(j=1,2,\cdots,M)$, it is always possible to choose $\DIS N_j^{(i)}$ such
that
\begin{equation}
\frac{\chi_{p_2}(j+1)}{\chi_{p_2}(j)} \gg
\frac{\chi_{p_1}(j+1)}{\chi_{p_1}(j)}.
\label{inequality}
\end{equation}
Then (\ref{4term}) is expanded as
\begin{eqnarray*}
\mbox{ (\ref{4term})}&=&\frac{(p_1-p_2)}{p_1p_2}\prod_{k=1}^2
\prod_{t'}^{t-1} \prod_{n'}^{n-1}
\left(1-\frac{p_k}{1+\delta_{t'}}\right)
\left(1-\frac{p_k}{1+\gamma_{n'}}\right)^{-1}
\\
&& \times \sum_{i=0}^{\infty}\sum_{i'=0}^{\infty} \left( a_{i,i'}
\chi_{p_1}(j+i)\chi_{p_2}(j+1+i') - b_{i,i'}
\chi_{p_2}(j+i)\chi_{p_1}(j+1+i')
\right),
\end{eqnarray*}
where the coefficients $\DIS a_{i,i'}$ are defined by 
 \begin{eqnarray*}
&&\left( \frac{p_1 p_2^2}{(p_1-q_{\ell_2}^{(2)})(p_2-q_{\ell_1}^{(1)})}\right)
 \prod_{k=1,2}
 \prod_{t'=-T_0}^{t-1}
\left(1-\frac{q_{\ell_k}^{(k)}}{1+\delta_{t'}}\right)^{-1}
\prod_{n'=n}^{N_0}\left(1-\frac{q_{\ell_k}^{(k)}}{1+\gamma_{n'}}\right)^{-1} \\
&&\qquad \qquad =\sum_{i=0}^{\infty}\sum_{i'=0}^{\infty} a_{i,i'}\left(
\frac{q_{\ell_1}^{(1)}}{p_{\ell_1}}
\right)^i
\left(  \frac{q_{\ell_2}^{(2)}}{p_{\ell_2}} \right)^{i'},
\end{eqnarray*}
 and $\DIS b_{i,i'} = \left(\frac{p_1}{p_2}\right) a_{i',i}$.
{}From (\ref{chiIEQ}), we evaluate
\begin{eqnarray*}
a_{0,0}\chi_{p_1}(j)\chi_{p_2}(j+1) &\gg &\sum_{i=0}^{\infty}
\sum_{\stackrel{\scriptstyle i'=0}{ i+i'\ne 0}}^{\infty}
a_{i,i'}\chi_{p_1}(j+i)\chi_{p_2}(j+1+i') \\
b_{0,0}\chi_{p_2}(j)\chi_{p_1}(j+1) &\gg &\sum_{i=0}^{\infty}
\sum_{\stackrel{\scriptstyle i'=0}{ i+i' \ne 0}}^{\infty}b_{i,i'}
\chi_{p_2}(j+i)\chi_{p_1}(j+1+i').
\end{eqnarray*}
Then, noticing $\DIS a_{0,0}=p_2$, $\DIS b_{0,0}=p_1$ and using
(\ref{inequality}), we find
\begin{eqnarray}
&&\lim_{\ve \to +0} \ve \log \tau(t+1,n+1,j+1) =
\max\left[0,K^{(1)}(t,n,j), K^{(2)}(t,n,j), \right. \nonumber \\
&& \qquad \qquad \qquad \left. K^{(1)}(t,n,j)+K^{(2)}(t,n,j)-A(j)\right],
 \\
&&K^{(i)}(t,n,j) := K_0^{(i)}-\sum_{j'=1}^j \ell_{j'}^{(i)}-\sum_{t'}^t
\min\left[\kappa_{t'},L^{(i)}\right] \nonumber \\
&& \qquad \qquad \qquad +\sum_{n'}^n \min\left[\theta_{n'},L^{(i)} \right]
\;\;(i=1,2),\\
&&A(j):=  L^{(2)}+\ell_{j+1}^{(2)} \;\; (0\le j \le M-1). 
\end{eqnarray}
This gives a two soliton solution.
For the scattering where
the larger soliton overtakes the smaller one like (I) in Theorem \ref{th:mainhomo},
the integer $\DIS \ell_j^{(1)}$ $(1 \le j\le M)$ corresponds
to the number of $j$ th balls in the larger soliton at
$t \to -\infty$, and $\DIS \ell_j^{(2)}$ corresponds to that of the smaller
soliton at $t \to +\infty$.
Since the balls in the smaller soliton at $t \to +\infty$
must be included in the larger soliton at $t \to -\infty$, the
condition $\DIS \ell_j^{(1)} \ge \ell_j^{(2)}$ must hold for soliton solutions.
Similarly, for the scattering where
the smaller soliton overtakes the larger one like (II) in Theorem \ref{th:mainhomo},
the integer $\DIS \ell_j^{(2)}$ $(1 \le j\le M)$ corresponds
to the number of $j$ th balls in the smaller soliton at
$t \to -\infty$, and $\DIS \ell_j^{(1)}$ corresponds to that of the larger
soliton at $t \to +\infty$.
We should also note that there are several freedoms to choose the
`phase' $A(j)$ in taking the ultradiscrete limit.
However we conjecture that
the above choice will cover all the canonical systems,
hence essentially all the time development patterns for $N=2$.

The $N$ soliton solution (\ref{eq:result}) is obtained in the same way.
The key in the construction is to evaluate the expansion:
\begin{eqnarray*}
&&\langle vac| \psi(p_1)\psi^*(q_1)\psi(p_2)\psi^*(q_2)\cdots
\psi(p_r)\psi^*(q_r) |vac \rangle \\
&&=\langle vac| \psi(p_1)\psi(p_2)\cdots
\psi(p_r)\psi^*(q_r)\psi^*(q_{r-1})\cdots \psi^*(q_r) |vac \rangle \\
&&=\frac{\prod_{1 \le i<j\le r}(p_i-p_j)(q_j-q_i)}{\prod_{1 \le i,j \le
r}(p_i-q_j)} \\
&&=\frac{\prod_{1 \le i<j\le r}(p_i-p_j)}{\prod_{i=1}^r p_i^r} \left(
q_r^{r-1} q_{r-1}^{r-2} \cdots q_2 + \mbox{\rm other terms}
\right)
\end{eqnarray*}
and show that this term gives the phase factor $A(\vec{\mu};j)$ and the
``other terms" do not contribute to the final results.
This can be done in the same manner as in the case of two soliton solutions.
We take
\begin{equation}
g(\t)=\prod_{i=1}^N \left(1+\psi(p_i,\t)\phi^*(p_i,\t)\right),
\end{equation}
where
\begin{equation}
\psi(p_i)=\psi(L^{(i)}:\ve), \quad \phi^*(p_i)=\phi^*(L^{(i)};
\{\ell_j^{(i)}\}; K_0^{(i)}: \ve) \quad (i=1,2,\cdots,N).
\end{equation}
We suppose
$$
L^{(1)} \ge L^{(2)} \ge \cdots \ge L^{(N)},
$$
and
$$
\ell_j^{(1)} \ge \ell_j^{(2)} \ge \cdots \ge \ell_j^{(N)}, \quad
(j=1,2,\cdots, M).
$$
Note that this implies: $\DIS    p_1>p_2>\cdots >p_N.$
The latter condition is also a natural constraint for $N$ soliton solutions
as in the case of two soliton solutions.
Finally we find that the result is  given by
(\ref{eq:result}).

\end{document}